\documentclass[12pt]{amsart}

\makeatletter
\def\section{\renewcommand{\@secnumfont}{\bf}%
\@startsection{section}{1}%
  \z@{.7\linespacing\@plus\linespacing}{.5\linespacing}%
  {\normalfont\larger[2]\bfseries\centering}}
\makeatother

\usepackage[backref=page]{hyperref}
\usepackage{color}

  \usepackage[nobysame]{amsrefs}
  \usepackage{amsmath,amssymb}

    \newcounter{talk}
    \newtheorem{theorem}{Theorem}[talk]
    \newtheorem{definition}[theorem]{Definition}
    \newtheorem{proposition}[theorem]{Proposition}
    \newtheorem{conjecture}[theorem]{Conjecture}
    \newtheorem{problem}[theorem]{Problem}
    \newtheorem{question}[theorem]{Question}
    \newtheorem{remark}[theorem]{Remark}
    \newtheorem{corollary}[theorem]{Corollary}

    \renewcommand{\paragraph}[1]{\bigbreak\noindent\textbf{#1}}

  \newcommand{\headline}[3]{%
    \par\bigbreak\bigskip
    \refstepcounter{talk}%
 \pdfbookmark[2]{#3}{talk.\thetalk}
    \noindent
  \addtocontents{toc}{\protect\contentsline {subsection}{\protect\tocsubsection{}{}{#1}}{\thepage}{talk.\thetalk}}%
  \addtocontents{toc}{\protect\talktitle{#2}{talk.\thetalk}}%
     \vbox{
      \begin{flushleft}
        \textbf{#1:}\\[1mm]
        {\large\bf #2}
      \end{flushleft}
    }\\
    \nopagebreak
}
  
  \newcommand{\course}[2]{\headline{#2}{#1}}
  \newcommand{\recordedcourse}[2]{\headline{#2}{#1~\protect\recorded}}
  \newcommand{\talk}[2]{\headline{#2}{#1}}
   \newcommand{\recordedtalk}[2]{\headline{#2}{#1~\protect\recorded}}
  
  \newcommand{\pers}[2]{#1 (#2)}

  \newcommand{\CCC}{\mathbb{C}}
  \newcommand{\FFF}{\mathbb{F}}
  \newcommand{\NNN}{\mathbb{N}}
  \newcommand{\QQQ}{\mathbb{Q}}
  \newcommand{\RRR}{\mathbb{R}}
  \newcommand{\ZZZ}{\mathbb{Z}}
  \newcommand{\monorightarrow}{\rightarrowtail}
  \newcommand{\notion}[1]{\emph{#1}}

  \DeclareFontFamily{OT1}{rsfs}{}
  \DeclareFontShape{OT1}{rsfs}{n}{it}{<-> rsfs10}{}
  \DeclareMathAlphabet{\mathscr}{OT1}{rsfs}{n}{it}

\makeatletter
\renewcommand{\tocsection}[3]{%
  \indentlabel{\@ifnotempty{#2}{\ignorespaces#1 \bf #2.\ }}{\bf #3}}
\def\l@section{\@tocline{1}{4pt}{0em}{5pc}{}} 
\def\l@subsection{\@tocline{2}{4pt}{1em}{4em}{}} 
\makeatother
\newcommand{\talktitle}[2]{\hskip1em%
	\pdfstartlink attr {/C [0.9 0 0]} goto name {#2}
	\textit{#1}%
	\pdfendlink}

\def\boxit#1{\setbox0=
 \vbox{\kern3pt\hbox{\kern3pt
 #1\kern3pt}\kern3pt}
  \dimen0=\ht0
 \divide \dimen0 by 2
 \hbox{\lower\dimen0\vbox{\hrule\hbox to \wd0{\vrule\hss
      \vbox{\copy0}\hss\vrule}\hrule}}}

\newcommand{\playsymbol}{\href
{http://www.birs.ca/events/2013/5-day-workshops/13w5019/videos/}%
{\raise3pt\hbox{\smaller[2]\color{blue}\boxit{$\triangleright$}}}%
}
\newcommand{\recorded}{\playsymbol}

\begin{document}

  \begin{center}
    \begingroup
    \huge\bf
      Arithmetic Groups
    \endgroup
      \\[\smallskipamount]
    \begingroup\large\bf
      (Banff, Alberta, April 14--19, 2013)
    \endgroup
  
    \bigskip
    
    edited by
            Kai-Uwe Bux,
            Dave Witte Morris,
            Gopal Prasad,
            and
            Andrei Rapinchuk 
              
    \bigskip
  
  \end{center}
  
  \vskip2cm

\noindent 
\textbf{Abstract.} 
We present detailed summaries of the talks that were given
    during a week-long workshop on \emph{Arithmetic Groups} at the Banff International Research Station in April 2013,
     organized by 
     Kai-Uwe Bux (University of Bielefeld),
        Dave Witte Morris (University of Lethbridge),
        Gopal Prasad (University of Michigan),
        and
        Andrei Rapinchuk (University of Virginia).
The vast majority of these reports are based on abstracts that were kindly provided by the speakers.
Video recordings of lectures marked with the symbol\,\recorded\ are available online at
\begin{center}
\url{http://www.birs.ca/events/2013/5-day-workshops/13w5019/videos/}
\end{center}

\vbox{\vskip 2cm}

\tableofcontents

\newpage
  \section{Introduction and Overview}
  	\thispagestyle{plain}
	\markboth{\textsc{Arithmetic Groups}}{\textsc{Introduction and Overview}}
\begingroup
\noindent
    The theory of arithmetic groups deals with  groups of matrices whose
    entries are integers, or more generally, $S$-integers in a global
    field. This notion has a long history, going back to the work of
    Gauss on integral quadratic forms. The modern theory of arithmetic
    groups retains its close connection to number theory (for example,
    through the theory of automorphic forms) but also relies on a
    variety of methods from the theory of algebraic groups, particularly
    over local and global fields (this area is often referred to as the
    arithmetic theory of algebraic groups), Lie groups, algebraic
    geometry, and various aspects of group theory (primarily, homological
    methods and the theory of profinite groups). At the same time,
    results about arithmetic groups have numerous applications in
    differential and hyperbolic geometry (as the fundamental groups of
    many important manifolds often turn out to be arithmetic),
    combinatorics (expander graphs), and other areas. There are also
    intriguing  connections and parallels (which are currently not so
    well-understood) between arithmetic groups and other important
    classes of groups, such as Kac-Moody groups, automorphism groups of
    free groups, and mapping class groups.

    The objective of the workshop was to survey the most significant
    results in the theory of arithmetic groups obtained primarily in the
    last five years, in order to make the new concepts and methods
    accessible to a broader group of mathematicians whose interests are
    closely related to arithmetic groups. The workshop brought together
    34 mathematicians, from the world's leading experts to recent PhD
    recipients and graduate students, working on a variety of problems
    involving arithmetic groups. This resulted in very active exchanges
    between and after the lectures. The scientific program of the
    workshop consisted of 3 mini-courses (two 45-min lectures each), 17
    survey and research talks (30 or 45 minutes), one of which (by Bertrand
    Remy) was not planned in advance, a Q\&A session, and an open
    problem session. 
    
    The subjects of the mini-courses were: 
    {\it Pseudo-reductive groups and
    their arithmetic applications} (Brian Conrad), {\it Towards an
    arithmetic Kac-Moody theory} (Ralf K\"ohl), and {\it
    Homological finiteness properties of arithmetic groups in positive
    characteristic} (Kevin Wortman). The
    mini-course on the pseudo-reductive groups focused on arithmetic
    applications of the theory of pseudo-reductive groups, developed by
    Conrad, Gabber and Prasad, which include the proof of fundamental
    finiteness theorems (finiteness of the class number, finiteness of
    the Tate-Schafarevich set, etc.)\ for {\it all} algebraic groups over
    the fields of positive characteristic, not just reductive ones. The
    course on Kac-Moody groups contained a series of results extending
    the known properties of (higher rank) arithmetic groups, such as
    property $(T)$ and (super)rigidity, to Kac-Moody groups over rings.
The course on
    homological finiteness properties contained an account of a major
    breakthrough in the area --- the proof of the Rank Conjecture. 

    \medbreak

 Most major areas within the theory of arithmetic groups were
    represented by at least one talk at the workshop.
    Topics that were not discussed, but  should be included in
    the program of  future meetings on the subject, include:
    	\begin{itemize} \itemsep=\smallskipamount
	\item connections between the
    cohomology of arithmetic groups and the theory of automorphic forms,
    \item   the virtual positivity of the first Betti number of certain rank one
    lattices and the growth of higher Betti numbers, 
    \item the analysis on
    homogeneous spaces modulo arithmetic groups.
    \end{itemize}

We conclude this introduction with a brief overview of the main themes that were discussed in the lectures.
More details of the three mini-courses are in Section~\ref{MiniSect}, and the other lectures are described individually in Section~\ref{LectSect}. 

  \subsection*{(a) Structural and homological properties}
    In addition to the mini-course on homological finiteness properties
    for $S$-arithmetic subgroups of semi-simple algebraic groups in
    positive characteristic, which contained an exposition of the work of
    Bux-K\"ohl-Witzel-Wortman on the rank conjecture, there was a talk
    by Stefan Witzel on finiteness properties for proper actions of
    arithmetic groups.

    The talk of Ted Chinburg demonstrated how to use the Lefschetz
    Theorem from algebraic geometry to show, in certain situations, that a
    ``large" arithmetic group can be generated by smaller arithmetic
    subgroups.

    In his talk, Vincent Emery showed how to bound the torsion
    homology of non-uniform arithmetic lattices in characteristic zero.

  \subsection*{(b) Profinite techniques for arithmetic groups}
    The method of analyzing arithmetic groups via the study of their
    finite quotients has a long history. One aspect of this approach,
    known as the congruence subgroup problem, is focused on
    understanding the difference between the profinite completion
    and the congruence
    completion. This difference is measured by the congruence kernel. In his
    talk, Andrei Rapinchuk gave a survey of the concepts and results
    pertaining to the congruence subgroup problem.

    In their talks,  which together virtually constituted another
    mini-course, Benjamin Klopsch and Christopher Voll presented their
    new results on the representation growth of $S$-arithmetic groups
    satisfying the congruence subgroup property (i.e., for which the
    congruence kernel is finite). These results are formulated in terms
    of the representation zeta function of the $S$-arithmetic group, and
    show, in particular, that the abscissa of
    convergence depends only on the root system.

    The talk of Pavel Zalesskii was devoted to the general question of
    when the profinite topology on an arithmetic group should be
    considered to be ``strong." The cohomological aspect of this
    question boils down to the notion of ``goodness" introduced by
    Serre. The central conjecture here asserts that if an $S$-arithmetic
    subgroup in characteristic zero fails the congruence subgroup
    property, then it  should be good, and the talk contained an account
    of the results supporting this conjecture.

  \subsection*{(c) Connections with Kac-Moody, automorphism groups of free groups, and  mapping class groups}
    As we already mentioned, Ralf K\"ohl gave a mini-course on Kac-Moody
    groups, in which he described how Kac-Moody groups are obtained from
    Chevalley groups by amalgamation, defined the Kac-Peterson topology
    on them, and established their important properties, including
    property $(T)$ and a variant of superrigidity. 
    The lectures generated so much interest in Kac-Moody groups among the
    conference participants that Bertrand R\'emy was asked to give a survey talk on the subject. He discussed various approaches to Kac-Moody groups, and
    stressed the utility of buildings in their analysis.

    In his talk, Alan Reid showed how methods from Topological Quantum 
    Field Theory can be used to prove that every finite group is a
    quotient of a suitable finite index subgroup of the mapping class
    group (for any genus).

    Lizhen Ji reported on the construction of complete geodesic metrics
    on the outer space~$X_n$ that are invariant under $\mathop{\mathrm{Out}}(F_n)$, where
    $F_n$ is the free group of rank~$n$.

  \subsection*{(d) Applications to geometry, topology, and beyond}
    Mikhail Belolipetsky gave a survey of the long line of research
    on hyperbolic reflection groups. (These are discrete isometry
    groups that are generated by reflections of hyperbolic $n$-space). 
    One of the
    central results here is that there are only finitely many conjugacy
    classes of maximal arithmetic hyperbolic reflection groups. This
    opens the possibility of classifying such groups.

    Matthew Stover reported on his results that provide a lower bound on 
    the number of
    ends (cusps) of an arithmetically defined hyperbolic manifold/orbifold
    or a locally symmetric space. In particular,
    one-end arithmetically defined hyperbolic $n$-orbifolds do not exist
    for any $n > 31$. The question about the existence of one-ended
    nonarithmetic finite volume hyperbolic manifolds remains wide open.

    T.~N.~Venkataramana spoke about the monodromy groups associated with
    hypergeometric functions. One of the central questions is when the
    monodromy group is a finite index subgroup of the corresponding
    integral symplectic group. A criterion was described for this to be
    the case.

  \subsection*{(e) Rigidity}
    Generally speaking, a group-theoretic rigidity theorem asserts that,
    in certain situations, any abstract homomorphism of a special
    subgroup (e.g, an arithmetic subgroup or a lattice) of the group of
    rational points of an algebraic group (resp., a Lie group, or a
    Kac-Moody group) can be extended to an algebraic (resp., analytic or
    continuous) homomorphism of the ambient group. The pioneering and most unexpected
    result of this type is Margulis's Superridity Theorem for
    irreducible lattices in  higher rank semi-simple Lie groups. As was
    pointed out by Bass, Milnor, and Serre, rigidity statements can also
    be proved for $S$-arithmetic groups if the corresponding congruence
    kernel is finite. 
    
    In his mini-course, Ralf K\"ohl showed how to
    use the known rigidity results for higher rank arithmetic groups to
    prove a rigidity theorem for Kac-Moody groups over~$\ZZZ$.

    In his talk, Igor Rapinchuk discussed his rigidity results for the
    finite-dimensional representations of elementary subgroups of
    Chevalley group of rank $> 1$ over arbitrary commutative rings.
    These results settle a conjecture of Borel and Tits about abstract
    homomorphisms for split algebraic groups of rank $> 1$ over fields
    of characteristic $\neq 2, 3$, and also have applications to
    character varieties of some finitely generated groups.

  \subsection*{(f) Weakly commensurable groups and connections to algebraic groups}
    The notion of weak commensurability for Zariski-dense subgroups of
    the group of rational points of semi-simple algebraic groups over a
    field of characteristic zero was introduced by G.~Prasad and
    A.~Rapinchuk. They were able to provide an almost
    complete answer to the question of when two $S$-arithmetic subgroups
    of absolutely almost simple algebraic groups are weakly
    commensurable. Using Schanuel's conjecture from transcendental
    number theory, they connected this work to the analysis of
    length-commensurable and isospectral locally symmetric spaces, and
    in fact obtained new important  results about isospectral spaces. In
    his talk, Rajan reported on his work which is based on a new notion
    of representation equivalence of lattices. While the condition of
    representation equivalence is generally stronger than
    isospectrality, it enables one to obtain results about
    representation equivalent locally symmetric spaces without using
    Schanuel's conjecture (which is still unproven).

    The work of Prasad-Rapinchuk also attracted attention to a wide
    range of questions in the theory of algebraic groups asking about a
    possible relationship between two absolutely almost simple algebraic
    groups $G_1$ and $G_2$ over the same field $K$, given the fact that
    they have the same isomorphism/isogeny classes of maximal $K$-tori.
    In his talk, Vladimir Chernousov reported on the recent results on a
    related problem of characterizing finite-dimensional central
    division algebras over the same field $K$ that have the same
    isomorphism classes of maximal subfields. He also formulated a
    conjecture that would generalize these results to arbitrary
    absolutely almost simple groups, and indicated that the results on
    division algebras enable one to prove this conjecture for inner
    forms of type $\textsf{A}_n$.

  \subsection*{(g) Applications to combinatorics}
    In his talk, Alireza Salehi Golsefidy discussed applications of
    arithmetic groups and their Zariski-dense subgroups to the
    construction of highly connected but sparse graphs known as {\it
    expanders}. It was pointed out by G.~A.~Margulis that families of
    expanders can be constructed from a discrete Kazhdan group $\Gamma$
    by fixing a finite system $S$ of generators of $\Gamma$ and
    considering the Cayley graphs $\mathrm{Cay}(\Gamma/N , S)$ of the finite
    quotients of $\Gamma$ with respect to this generating system. Later,
    using deep number-theoretic results, Lubotzky, Phillips, and Sarnak
    showed that one also obtains a family of expanders from the
    non-Kazhdan group $SL_2(\ZZZ)$ by fixing a system $S$ of
    generators of the latter and considering the Cayley graphs
    $\mathrm{Cay} \bigl( SL_2(\ZZZ/d\ZZZ) , S \bigr)$ for the congruence quotients.

    Lubotzky raised the question of whether one gets a family of
    expanders if one takes, for example, an arbitrary finitely generated
    Zariski-dense subgroup $\Gamma \subset SL_2(\ZZZ)$, and
    considers the Cayley graphs of the congruence quotients
    $SL_2(\ZZZ/d\ZZZ)$ with respect to a fixed finite
    generating set of the subgroup. (It is known that $\Gamma$ will map
    surjectively onto $SL_2(\ZZZ/d\ZZZ)$ for all $d$ prime
    to some $d_0$ that depends on $\Gamma$.)

    Alireza Salehi Golsefidy surveyed the important progress on this question in the
    context of general arithmetic groups, including the groundbreaking work of
    Bourgain and Gamburd, and his recent results with Varj\'{u}.

\endgroup
  
  \newpage\addtocontents{toc}{\protect\vfil\protect\bigskip}
  \section{Mini-Courses} \label{MiniSect}
  	\thispagestyle{plain}
	\markright{\textsc{Mini-Courses}}

 \begingroup
\recordedcourse{Pseudo-reductive groups and their arithmetic applications}
       {\pers{Brian~Conrad}{Stanford University}}
       {Conrad: Pseudo-reductive groups and their arithmetic applications}
    The talk first described some motivation for, and examples of, the theory
    of pseudo-reductive groups \cite{Con-pred},
    culminating in the main structure
    theorem (ignoring subtleties that arise in characteristics 2 and 3).
    Then, it discussed how to apply the main structure theorem
    in the context of proving some finiteness theorems over
    global function fields \cite{Con-finite}.
   Such results include affirmative solutions to
    open questions which render unconditional several results in
    \cite{Con-oesterle} on the behavior of sizes of degree-1
    Tate--Shafarevich sets and Tamagawa numbers for linear algebraic
    groups over global function fields.

    A very highly-recommended survey of both the general theory and
    arithmetic applications is given in the Bourbaki report
    \cite{Con-remy} by Bertrand R\'emy. This includes a user-friendly
    overview of the contents of \cite{Con-pred}, indicating where main
    results can be found and how the logical development of the main
    proofs proceeds.

    \paragraph{Definition of pseudo-reductivity:}
    triviality of the so-called {\em $k$-unipotent radical}
    $\mathscr{R}_{u,k}(G)$ (the largest smooth connected normal
    unipotent $k$-subgroup of a smooth connected affine $k$-group), with
    $k$ a general field.  The formation of $\mathscr{R}_{u,k}(G)$
    commutes with separable extension on $k$, including $K \rightarrow
    K_v$ for a global function field $K$ and place $v$. We say $G$ is
    {\em pseudo-reductive} if it is smooth connected affine and
    $\mathscr{R}_{u,k}(G)=1$.    This coincides with ``connected
    reductive'' when $k$ is perfect. Several examples were given over
    any imperfect field, both commutative as well as perfect ($G =
    \mathscr{D}(G)$). The most basic example is the Weil restriction
    ${\rm{R}}_{k'/k}(G')$ for a finite extension $k'/k$ and a connected
    reductive $k'$-group $G'$; this is never reductive if $G' \ne 1$ and
    $k'$ is not separable over $k$.

    \paragraph{Surprises:} the purely inseparable Weil restriction
    may fail to preserve dimension (e.g.,
    ${\rm{R}}_{k'/k}(\mu_p)$ has positive dimension for a nontrivial
    purely inseparable finite extension $k'/k$ in characteristic $p >
    0$). It may also fail to preserve surjectivity or perfectness:
    ${\rm{R}}_{k'/k}({\rm{PGL}}_p)$ has a nontrivial commutative
    quotient modulo the image of ${\rm{R}}_{k'/k}({\rm{SL}}_p)$ for a
    nontrivial purely inseparable finite extension $k'/k$ in
    characteristic $p > 0$. Other surprises: pseudo-reductivity can be
    lost under central quotient, and also some quotients by smooth
    connected normal $k$-subgroups.

    \paragraph{Good news:} Cartan $k$-subgroups are always commutative,
    ${\rm{R}}_{k'/k}(G')$ is perfect when $G'$ is simply connected, and
    there is a theory of root systems (which can be non-reduced in
    characteristic 2, even if $k = k_s$). A related concept is
    pseudo-semisimplicity; this has some surprises
    (there are two possible definitions, one of which is ``wrong'').
    Overall, pseudo-reductivity is not a particularly robust
    concept, so its main purpose is the role it plays in trying to prove
    theorems about rather general linear algebraic groups when one might
    not have much control (e.g. for Zariski closures, working with
    stabilizer schemes for a group action, etc.).

    \paragraph{A general principle:} we cannot hope to understand the commutative
    pseudo-reductive groups, but we will aim to describe the general
    structure modulo that ignorance. Life could be worse; at least we
    are remaining ignorant only about something commutative.

    \paragraph{The standard construction:}
    this is a procedure which is a kind of pushout that replaces
    a Weil-restricted maximal torus ${\rm{R}}_{k'/k}(T')$ from a simply
    connected semisimple $k'$-group $G'$ with another commutative
    pseudo-reductive group $C$ according to a very specific kind of
    procedure.   The final output of this process is a central quotient
    presentation
    $$G = ({\rm{R}}_{k'/k}(G') \rtimes C)/{\rm{R}}_{k'/k}(T').$$
    This can be generalized to allow several extensions $k'_i/k$, and
    admits a precise uniqueness aspect as well, determines all of the
    data $(G', k'/k, T', C)$ in terms of a choice of maximal $k$-torus
    of $T$.

    \paragraph{A general principle} for
    applying the structure theory of pseudo-reductive groups: if a
    theorem is known in the smooth connected solvable affine case over
    $k$ and in the connected semisimple case over all finite extensions
    of $k$ then ``probably'' one can use the structure theorem via
    standard presentations (plus extra care in characteristics 2 and 3)
    to prove the result for all smooth connected linear algebraic groups
    (and something without smoothness or connectedness, depending on the
    specific assertion).

    \paragraph{Some of the finiteness questions} one would like to settle (all of
    which have long been known in the affirmative in the connected
    reductive case): finiteness for ``class numbers'' of smooth
    connected linear algebraic groups over global function fields,
    finiteness for degree-1 Tate--Shafarevich sets for all affine
    algebraic group schemes over global function fields, finiteness for
    Tamagawa numbers of smooth connected affine groups, and finiteness
    for obstruction sets to a local-global principle for orbits over
    global function fields (i.e., if $x, x' \in X(k)$ are in the same
    $G(k_v)$-orbit for all $v \not\in S$ then the images of $x$ and $x'$
    in $G(k)\backslash X(k)$ might not coincide but are at least
    constrained within a finite set, depending on $S$).  In this final
    question it is important that we do not impose smoothness hypotheses
    (or any hypotheses at all) on the stabilizer schemes for the
    geometrically transitive action.
 %
    The finiteness
    question for orbits reduces to finiteness for Tate--Shafarevich sets
    for the {\em stabilizer} scheme, so it is important for the latter
    finiteness result that we do not demand smoothness hypotheses
    (but they can be imposed by a trick).

    \paragraph{Example applications:} finiteness for Tate--Shafarevich
    sets can be reduced to the pseudo-reductive case, and the {\em form}
    of the ``standard presentation'' is very well-suited to
    then pulling up the known result for the semisimple and commutative
    cases, essentially using vanishing theorems for simply connected
    groups.   Second: finiteness for Tamagawa numbers is settled in
    general by a different kind of argument with the ``standard
    presentation'' for pseudo-reductive groups to pull it up from the
    known semisimple and commutative cases.
    Third: the original finiteness question for
    the local-global principle with orbits.  That indeed works out
    affirmatively, due to the established case for Tate--Shafarevich
    sets. Fourth: various formulas for the
    behavior of Tamagawa numbers under exact sequences that were proved
    in Oesterl\'e's paper \cite{Con-oesterle} conditionally on certain
    unknown finiteness results are now all valid unconditionally.
\endgroup

\begingroup
  \newcommand{\SL}{\operatorname{SL}}
\course{Kac-Moody groups}
       {\pers{Ralf~K{\"o}hl}{University of Gie{\ss}en}}
       {Koehl: Kac-Moody groups}
%
  By theorems of Tits and Curtis,
  a Chevalley group over a field with at least four elements
  is the product of its rank two subgroups amalgamated
  over the rank one subgroups. The following result goes to
  show that in the context of Chevalley groups over local fields,
  the topology is also forced by the rank one subgroups:

  \begin{theorem}[Gl{\"o}ckner-Hartnick-K{\"o}hl \cite{MR2684413}]
    Let $\FFF$ be a local field, and let $G$ be
    a Chevalley group over $\FFF$. Then the Lie group
    topology on $G$ is the finest group topology making the
    embeddings of the fundamental rank one subgroups (endowed with
    their Lie group topologies) continuous.
  \end{theorem}

  Now let $\Delta$ be a $2$-spherical Dynkin diagram without loops and
  let $\FFF$ be a field with at least four elements. For each node
  $\alpha\in\Delta$ let $G_\alpha$ be a copy of $\SL_2(\FFF)$;
  and for each pair $\alpha,\beta\in\Delta$ let $G_{\alpha,\beta}$
  be a simply connected Chevalley group over $\FFF$ of the type
  given in $\Delta$. There are obvious inclusions
  $G_\alpha \monorightarrow G_{\alpha,\beta}$. The Kac-Moody
  group $G_\Delta(\FFF)$ can then be described as the product
  of the $G_{\alpha,\beta}$ amalgamated over the $G_\alpha$ and
  is uniquely determined by $\Delta$ since $\Delta$ does not contain
  loops. Every $2$-spherical split Kac-Moody group arises this way.
  \begin{definition}
    The \notion{Kac-Peterson} topology on $G_\Delta(\FFF)$ is
    the finest group topology that makes the canonical embeddings
    $G_\alpha \monorightarrow G_\Delta(\FFF)$ continuous.
  \end{definition}
  \begin{theorem}[Hartnick-K{\"o}hl-Mars \cite{koehl-twin}]
    The group $G_\Delta(\FFF)$ with the Kac-Peterson topology is Hausdorff and a
    $k_\omega$-space, i.e., it is the direct limit of an ascending
    sequence of compact Hausdorff subspaces.
    If $\Delta$ is not spherical, then the Kac-Peterson topology
    is neither locally compact nor metrizable.
  \end{theorem}
  In particular, the existence of a Haar measure is not guaranteed
  for non-spherical Kac-Moody groups with the Kac-Peterson topology.
  \begin{theorem}[Hartnick-K{\"o}hl \cite{koehl-kazhdan}]
    Let $\FFF$ be a local field and let $G_\Delta$ be an irreducible (i.e., $\Delta$
    is connected),
    $2$-spherical split Kac-Moody group. Then $G_\Delta(\FFF)$ with the
    Kac-Peterson topology has Kazhdan's property~$(T)$.
  \end{theorem}
  The subgroup $G_\Delta(\ZZZ)$ is discrete and finitely generated. One would
  like to think of $G_\Delta(\ZZZ)$ in analogy to arithmetic lattices. It
  is, however, an open question whether it inherits property~$(T)$ from
  $G_\Delta(\RRR)$.

  The following follows easily from a theorem of Caprace and
  Monod on Chevalley groups acting on {\small CAT(0)} polyhedral
  complexes applied to the Davis realization of the twin building
  for $G_\Delta(\RRR)$:
  \begin{proposition}
    Let $L$ be an irreducible Chevalley group of rank at least two,
    let $G_\Delta(\RRR)$ be a Kac-Moody group, and let
    $\varphi \colon L(\ZZZ) \rightarrow G_\Delta(\RRR)$
    be a group homomorphism.
    Then the image $\varphi( L(\ZZZ) )$ is a bounded subgroup,
    i.e., it lies in the intersection of two parabolic subgroups
    of opposite sign. In particular, $\varphi( L(\ZZZ) )$
    is contained in an algebraic subgroup of $G_\Delta(\RRR)$.
  \end{proposition}
  This can be extended to yield an analogue of Margulis
  superrigidity:
  \begin{theorem}[Farahmand-Horn-K{\"o}hl]
    Let $G_\Delta(\RRR)$ and $G_{\Delta'}(\RRR)$ be irreducible $2$-spherical
    Kac-Moody groups and let
    $\varphi \colon G_\Delta(\ZZZ) \rightarrow G_{\Delta'}(\RRR)$
    be a group homomorphism with Zariski dense image. Then there
    exists $n\in\NNN$ such that the restriction of $\varphi$ to
    $G_\Delta(n\ZZZ)$ extends uniquely to a continuous homomorphism
    \(
      G_\Delta(\RRR) \rightarrow G_{\Delta'}(\RRR)
    \)
    with respect to the Kac-Peterson topologies.
  \end{theorem}
  The proof proceeds by first dealing with the case that $G_\Delta$ is
  a Chevalley group, where the Kac-Peterson topology is the
  Lie group topology. For general $G_\Delta$, the statement is reduced
  to the rank two subgroups, which are Chevalley groups. Using the
  the presentations of $G_\Delta$ and $G_{\Delta'}$ as products of their respective
  rank two subgroups amalgamated along their rank one subgroups,
  one constructs the extension of $\varphi$. Since the
  Kac-Peterson topology is universal with respect to the Lie topologies
  on the rank one subgroups, it follows that $\varphi$ is continuous.
\endgroup

\begingroup
\recordedcourse{Finiteness properties of arithmetic groups over function fields}
       {\pers{Kevin~Wortman}{University of Utah}}
       {Wortman: Finiteness properties of arithmetic groups over function fields}
  Recall that a group $\Gamma$ has
  finiteness length $\leq m$ if it has a classifying space whose
  $m$-skeleton is finite. In this case, the cohomology of $\Gamma$
  is clearly finitely generated in dimensions $\leq m$.

  Let $K$ be a global function field of characteristic $p > 0$,
  let $\FFF_p$ be the finite field with $p$~elements,
  and let $\mathbf{G}$ be a
  connected noncommutative absolutely almost simple $K$-isotropic
  $K$-group. Let
  \(
    d :=
    \sum_{p \in S} \operatorname{rk}_{K_p}(\mathbf{G})
  \)
  denote the sum of the local ranks of $\mathbf{G}$.
  With this notation fixed,
  the two main results are:
  \begin{theorem}[``Rank Theorem'', Bux-K\"ohl-Witzel \cite{MR2999042}]
      The finiteness length
      $\phi(\Gamma)$ of the $S$-arithmetic subgroup $\Gamma=
      \mathbf{G}(\mathcal{O}_S)$ is $d-1$.
  \end{theorem}
  \begin{theorem}[Wortman] \label{Wortman-NotFG}
      For some subgroup $\Gamma'$ of finite index in
      $\Gamma$, the cohomology ${\rm H}^d(\Gamma',\FFF_p)$ is not finitely
      generated.
  \end{theorem}
  
   (At the time of the conference, a mild restriction
      on the $K$-type of $\mathbf{G}$ was needed in Theorem~\ref{Wortman-NotFG}, but Wortman was soon able to remove this restriction.)
      
  Results on finiteness properties of arithmetic groups have a
  long history. The Euclidean algorithm shows that
  $\operatorname{SL}_n(\ZZZ)$
  is finitely generated. Finite presentability of these groups
  is a classical
  application of Siegel domains. Raghunathan~\cite{MR0230332} proved that arithmetic
  groups in characteristic~$0$ enjoy all finiteness properties. In fact,
  he showed that they have a torsion-free subgroup of finite index that
  is the fundamental group of a compact aspherical manifold with boundary.
  Borel-Serre \cite{MR0447474} have shown that $S$-arithmetic subgroups of reductive
  groups in characteristic~$0$ also enjoy all finiteness properties.

  The picture in positive characteristic is different. Nagao \cite{MR0114866} showed
  that $\operatorname{SL}_2(\FFF_q[t])$ is not even finitely
  generated. Behr~\cite{MR0238853} proved that $\Gamma$ as in the Rank Theorem
  is finitely generated if and only if $d>1$. Stuhler~\cite{MR568936} showed
  that $\operatorname{SL}_2(\mathcal{O}_S)$ has finiteness length
  $|S|-1=d-1$. Abels~\cite{MR1177335} and Abramenko~\cite{MR1286827}
  independently showed
  that $\operatorname{SL}_n(\FFF_q[t])$ has finiteness length
  $n-2=d-1$, provided $q$ is large enough. Sometime during the 1980s,
  the pattern became transparent. Behr turned it into a serious
  conjecture when he proved in \cite{MR1603845} that the $S$-arithmetic subgroup
  $\Gamma$ of the Rank Theorem is finitely presented if and only if $d>2$.

  A significant step toward the Rank Theorem was the proof of its
  ``negative half'' by Wortman and Bux in \cite{MR2270455}, where they showed
  that $\phi(\Gamma)<d$. In 2008 (published in \cite{MR2819164}),
  Wortman and Bux also settled the Rank Theorem in full for
  groups of global rank~$1$. The major improvement was a geometric
  filtration of the Bruhat-Tits building for $\Gamma$ defined by
  Busemann functions. The relative links of this filtration are
  larger than those occurring in combinatorially defined filtrations
  used previously: the new relative links
  are hemi-sphere complexes in spherical buildings, whose connectivity
  properties have been established by Schulz~\cite{MR3039769}. The proof of
  the Rank Theorem for arbitrary groups follows this line of thought.
  Here, Behr-Harder reduction theory is the source of the Busemann
  functions and the associate filtration.

  Theorem~\ref{Wortman-NotFG} above is a considerable strengthening of the negative half of the
  Rank Theorem.
  For
  \(
    \operatorname{SL}_2(\mathcal{O}_S)
  \),
  Stuhler succeeded in proving that the homology in the
  critical dimension (here $d = |S|$) is infinitely generated, by using of a spectral
  sequence argument. In the other works cited above, the finiteness
  length was deduced by combinatorial or geometrical means that do
  not detect homology in the critical dimension.

  The main difficulty is that the action of $\Gamma$ on its
  associated Bruhat-Tits building $X$ is not free; in fact, the
  orders of the point stabilizers are not bounded. Wortman uses the
  height function from the Rank Theorem to pass to a cocompact
  subspace $X(0)$, which is $(d-2)$-connected. Gluing in free
  $\Gamma$-orbits of cells of
  dimensions $d$ and $d+1$, he obtains a $d$-connected space $Y$
  on which $\Gamma$ acts with stabilizers of uniformly bounded
  order, and he obtains a $\Gamma$-equivariant map $Y \rightarrow X$.
  He can now pass to a finite index subgroup $\Gamma'\leq\Gamma$
  that acts freely on $Y$. He then constructs an infinite family
  of cocycles on $Y$ and a ``dual'' family of cycles on $X$ paired
  via the comparison map $Y \rightarrow X$.
  The supports of the cycles in $X$ increase in height and
  ``low'' cocycles evaluate trivially on higher cycles. On the
  other hand, each cocycle evaluates non-trivially on the
  corresponding cycle. This shows that the cocycles are non-trivial
  and linearly independent.
\endgroup

  \newpage\addtocontents{toc}{\protect\newpage}
\section{Research Lectures and Survey Talks} \label{LectSect}
  	\thispagestyle{plain}
	\markright{\textsc{Research Lectures and Survey Talks}}

\begingroup
\recordedtalk{Arithmetic hyperbolic reflection groups}
     {\pers{Mikhail~Belolipetsky}{IMPA, Brazil}}
     {Belolipetsky: Arithmetic hyperbolic reflection groups}
  A \emph{hyperbolic reflection group} $\Gamma$ is a
  discrete subgroup of the group of isometries of the hyperbolic
  $n$-space $\operatorname{Isom}(\mathbf{H}^n)$ generated by
  reflections in the faces of a hyperbolic polyhedron $P\subset
  \mathbf{H}^n$. If $P$ has finite volume, then $\mathbf{H}^n/\Gamma =
  \mathcal{O}$ is a finite volume hyperbolic orbifold, which is
  obtained by ``mirroring'' the faces of $P$. A reflection group
  $\Gamma$ is \emph{maximal} if there does not exist a reflection
  group $\Gamma' \subset \operatorname{Isom}(\mathbf{H}^n)$ that
  properly contains $\Gamma$. A reflection group is called
  \emph{arithmetic} if it is an arithmetic subgroup of
  $\operatorname{Isom}(\mathbf{H}^n) = \mathrm{PO}(n,1)$.

  This talk was about finiteness results for maximal arithmetic
  hyperbolic reflection groups. After a brief review of previous
  foundational work by Vinberg and Nikulin, it focused attention on the
  results obtained in this area in the last 10 years. The important
  breakthrough was achieved by Maclachlan--Long--Reid \cite{Bel-LMR}
  and Agol \cite{Bel-A06}, who proved that there exist only finitely
  many conjugacy classes of arithmetic maximal hyperbolic reflection
  groups in dimensions $n=2$ and $n=3$, respectively. Later, this
  was proved for all dimensions:
  
  \begin{theorem}[Agol-Storm-Belolipetsky-Whyte \cite{Bel-ABSW}, Nikulin
  \cite{Bel-N07}]
  There are only finitely
  many conjugacy classes of arithmetic maximal hyperbolic reflection
  groups in any fixed dimension~$n$.
  \end{theorem}

  The finiteness theorem allows us to talk about a classification of 
  arithmetic hyperbolic reflection groups. Two types of problems are
  considered here: proving quantitative bounds for the invariants of
  the groups, and constructing examples that fit into the bounds. 
  Results in these directions were obtained in the
  papers \cites{Bel-B09, Bel-B11, Bel-Mcleod, Bel-M11, Bel-N11,
  Bel-BLi, Bel-BM}. The end of the talk discussed the main open
  problems in the area.
\endgroup

  \begingroup
  \recordedtalk{A finiteness theorem for the genus}
     {\pers{Vladimir Chernousov}{University of Alberta}}
     {Chernousov: A finiteness theorem for the genus}
  Given a finite-dimensional central
  division algebra $D$ over a field $K$, the {\it genus}
  $\mathbf{gen}(D)$ is defined to be the set of isomorphism classes of central
  division $K$-algebras having the same (isomorphism classes of)
  maximal subfields as $D$. 
  One would like to have
  qualitative and quantitative results for $\mathbf{gen}(D)$ over
  arbitrary fields. 
  
  This general question is related to other
  interesting problems in division algebras, quadratic forms, Galois
  cohomology and even differential geometry (a question of this nature
  was raised in the paper \cite{Raj-PR} on length-commensurable and
  isospectral locally symmetric spaces). Since every division algebra
  is the union of its maximal subfields, questions about the genus can
  informally be thought of questions about ways to use the same
  subfields and construct a different algebra (so, in some sense,
  these are analogs of the Banach-Tarski paradox for division
  algebras).

  More precisely, we have the following questions.
  \begin{question}\label{cher-q-1}
    When is $|\mathbf{gen}(D)| = 1$?
  \end{question}

  We note that $|\mathbf{gen}(D)| = 1$ if and only if $D$ is uniquely
  determined by its maximal subfields. Since $D$ and its opposite
  algebra $D^{\mathrm{op}}$ have the same maximal subfields, an
  affirmative answer to Question~\ref{cher-q-1} is  possible only if
  ${\rm exp}(D)=2$.

  \begin{question}\label{cher-q-2}
    When is $|\mathbf{gen}(D)|<\infty$?
  \end{question}

  Various people including Garibaldi, Rost, Saltman, Shacher,
  Wadsworth, have described a method for constructing non-isomorphic
  quaternion algebras over very large fields (having infinite
  transcendence degree over the prime field) with the same quadratic
  subfields, which actually shows that the genus of a quaternion
  algebra over such a field can be infinite. This suggests that
  Question~\ref{cher-q-2} should be considered primarily over finitely
  generated fields.

  Work of Vladimir Chernousov with Andrei and Igor Rapinchuk 
  developed a general approach to proving the finiteness of the
  genus of a division  algebra, and estimating its size, based on an
  analysis of the unramified Brauer group with respect to an
  appropriate set of discrete valuations of $K$. This approach yields,
  in particular, the following two theorems.

  \begin{theorem}[``Stability Theorem'', Chernousov-Rapinchuk-Rapinchuk \cites{Cher-CRR1,Cher-CRR2}]\label{cher-t-1}
    Let $K$ be a field of characteristic not $2$. If
    $|\mathbf{gen}(D)| = 1$ for any central division $K$-algebra $D$
    of exponent $2$, then the same is true for any division algebra
    of exponent $2$ over the field of rational functions $K(x)$.
  \end{theorem}
  \begin{corollary}
    Let $k$ be either a finite field of
    characteristic not $2$ or a~number field, and let
    $K=k(x_1,\ldots,x_r)$ be a finitely generated purely transcendental
    extension of $k$. Then, for any central division $K$-algebra $D$ of
    exponent $2$, we have $|\mathbf{gen}(D)|=1$.
  \end{corollary}
  \begin{theorem}[Chernousov-Rapinchuk-Rapinchuk \cites{Cher-CRR1,Cher-CRR2}]\label{cher-t-2}
    Let $K$ be a finitely generated field. If $D$ is a central
    division $K$-algebra of exponent prime to
    ${\rm char} \: K$, then $\mathbf{gen}(D)$ is finite.
  \end{theorem}

  Furthermore, the authors proposed a generalization of the notion of
  genus to arbitrary absolutely almost simple algebraic $K$-groups,
  based on the isomorphism classes of maximal $K$-tori. (Possible
  variations of this notion can be based on the consideration of 
  isogeny classes and/or  some special classes of maximal $K$-tori,
  e.g., generic tori.)

  In view of Theorem~\ref{cher-t-2}, the
  following seems natural.

  \begin{conjecture}[Chernousov-Rapinchuk-Rapinchuk \cites{Cher-CRR1,Cher-CRR2}] \label{cher-c}
    Let $G$ be an absolutely almost simple simply connected algebraic
    group over a~finitely generated field $K$ of characteristic zero
    (or of ``good" characteristic relative to $G$). Then there exists
    a \emph{finite} collection $G_1, \ldots , G_r$ of $K$-forms of
    $G$ such that if $H$ is a $K$-form of $G$ having the same
    isomorphism classes of maximal $K$-tori as $G$, then $H$ is
    $K$-isomorphic to one of the $G_i$'s.
  \end{conjecture}

  The proof of Theorem~\ref{cher-t-2} yields a
  proof of this conjecture for inner forms of type $\textsf{A}_{\ell}$.
\endgroup

\begingroup
\recordedtalk{Generating arithmetic groups by small subgroups using Lefschetz Theorems}
     {\pers{Ted~Chinburg}{University of Pennsylvania}}
     {Chinburg: Generating arithmetic groups by small subgroups using Lefschetz Theorems}
  The following result provides interesting examples of large arithmetic groups that are generated by specific small arithmetic subgroups.

  \begin{theorem}[Chinburg-Stover \cite{Chin-CS}] \label{thm:main}
    Let $G$ be
    \(
      \mathrm{SL}_2(\RRR) \times \mathrm{SL}_2(\RRR)
    \)
    or $\mathrm{SU}(2, 1)$ with corresponding hermitian symmetric
    space $X$, and let $\Gamma$ be a cocompact arithmetic lattice
    in $G$. Assume that the complex algebraic surface
    $S = \Gamma \backslash X$ contains a holomorphically immersed
    totally geodesic projective algebraic curve. Then there exist
    finitely many such curves $C_1, \dots, C_r \subset S$ and
    positive integers $\alpha_1, \dots, \alpha_r$ such that:
    \begin{enumerate}
      \item\label{Chin-a}
        The divisor $D = \sum \alpha_j C_j$ is a connected
        effective divisor on $S$.
      \item\label{Chin-b}
        The image of $\pi_1(|D|)$ in $\Gamma$ under the natural
        homomorphism is a finite index subgroup.
    \end{enumerate}
  \end{theorem}

  Statement~\ref{Chin-a}  is shown using the commensurator of
  $\Gamma$ and the fact that there are infinitely many commensurability
  classes of Fuchsian curves on $S$ once one exists.   The proof of
  statement~\ref{Chin-b}  uses work of Nori  \cite{Chin-N} and of Napier
  and Ramachandran \cite{Chin-NR} on Lefschetz Theorems
  for sufficiently positive divisors on complex varieties.

  One consequence of the theorem is the following structure theorem for
  the Albanese varieties of arithmetic complex hyperbolic $2$-manifolds.

  \begin{theorem}[Chinburg-Stover \cite{Chin-CS}] \label{thm:albanese}
    Let $\Gamma$ be as in Theorem \ref{thm:main} and suppose
    $G = \mathrm{SU}(2,1)$. There exists $r \geq 1$ and Fuchsian
    curves $C_1, \dots, C_r$ on
    $S = \Gamma \backslash \mathbf{H}_{\CCC}^2$ such that if
    the Albanese variety $\mathrm{Alb}(S)$ is nontrivial, then every
    simple factor of $\mathrm{Alb}(S)$ is isogenous to a factor of
    the Jacobian of the normalization $C_j^\#$ of (at least) one
    of the curves $C_j$.
  \end{theorem}

  If $\Gamma$ is a congruence arithmetic lattice of simple type,
  Gelbart and Rogawski \cite{Chin-GR} proved that the
  first cohomology group of $\Gamma \backslash \mathbf{H}_{\CCC}^2$,
  which determines
  $\mathrm{Alb}(\Gamma \backslash \mathbf{H}_{\CCC}^2)$, arises
  from the theta correspondence. Murty and Ramakrishnan then used
  Faltings' work on the Mordell conjecture
  in \cite{Chin-MR} to show that the simple factors
  of $\mathrm{Alb}(\Gamma \backslash \mathbf{H}_{\CCC}^2)$
  are, in fact, CM abelian varieties. This gave a positive answer to
  a question of Langlands.

  Theorem~\ref{thm:albanese} provides information of a different
  nature about
  $\mathrm{Alb}(\Gamma \backslash \mathbf{H}_{\CCC}^2)$ for both
  congruence and noncongruence $\Gamma$, namely that
  $\mathrm{Alb}(\Gamma \backslash \mathbf{H}_{\CCC}^2)$ is
  built from the Jacobians of Fuchsian curves on $\Gamma \backslash X$.
  By work of Kazhdan in \cite{Chin-K}, one can always find elements
  in the commensurability class with nontrivial Albanese variety,
  and there are always noncongruence groups in the commensurability
  classes under consideration. It would be interesting to know
  whether or not the factors which appear for noncongruence
  $\Gamma$ must also have complex multiplication.
\endgroup

\begingroup
\recordedtalk{Bounds for torsion homology of arithmetic groups}
     {\pers{Vincent~Emery}{Stanford University}}
     {Emery: Bounds for torsion homology of arithmetic groups}
  Let $\mathbf{G}$ be a connected semisimple real algebraic
  group, such that $\mathbf{G}(\RRR)$ has trivial center and no compact
  factor. A result  of Gelander
  \cite{Em-Gel04} can be used to prove the following theorem that bounds
  the torsion homology of nonuniform arithmetic lattices
  $\Gamma \subset \mathbf{G}(\RRR)$, without requiring  $\Gamma$ to be 
  torsion-free (but with a restriction on~$\mathbf{G}$). 
  Results of this type for Betti numbers hold in much greater
  generality (see, for instance, \cite{Em-Samet}).

  \begin{theorem}[Emery \cite{Em-EmeK2}]
    Let $\mathbf{G}$
    be as above, and assume, for all irreducible lattices $\Gamma \subset
    \mathbf{G}(\RRR)$, that we have $H_q(\Gamma,\QQQ) = 0$ for $q=1,\dots,j$. Then
    there exists a constant $C_{\mathbf{G}} > 0$, such that, for each
    irreducible nonuniform arithmetic lattice $\Gamma \subset
    \mathbf{G}(\RRR)$, the following bound on torsion homology holds:
    \[
    \log \left| H_j(\Gamma, \ZZZ) \right| \le
        C_{\mathbf{G}} \; \mathrm{vol} \bigl( \Gamma \backslash \mathbf{G}(\RRR) \bigr).
    \]
  \end{theorem}

 The following theorem about $K$-theory
  of number fields is obtained by combining Theorem~1
   with recent results of Calegari and Venkatesh
  \cite{Em-CalVen}. It improves --- for totally imaginary
  fields --- previous bounds due to Soul\'e \cite{Em-Soule03}.

  \begin{theorem}[Emery \cite{Em-EmeK2}]
    Let $d > 5$ be an
    integer. There exists a constant $C = C(d) > 0$, such that, for each
    totally imaginary field $F$ of degree~$d$,  we have:
    $$
        \log |K_2(\mathcal{O}_F) \otimes R| \le C |D_F|^{2} (\log|D_F|)^{d-1},
    $$
    where $R = \ZZZ[\frac{1}{6w_F}]$, $D_F$ is the discriminant, and $w_F$
    the number of roots of unity in $F$.
  \end{theorem}
\endgroup

\begingroup
\talk{Expansion properties of linear groups}
     {\pers{Alireza Salehi Golsefidy}{University of California, San Diego}}
     {Golsefidy: Expansion properties of linear groups}
%
  Highly connected sparse graphs are extremely useful in the theory of communication,
  theoretical computer science, and pure mathematics (see the beautiful
  surveys \cite{Gol-HLW}, \cite{Gol-Lub} and \cite{Gol-Kow}). A
  family $\{X_i\}_{i=1}^{\infty}$ of finite $k$-regular graphs is
  called a family  of {\em expanders} if, for some positive number~$c$, we
  have 
  	\[ \text{$c\min\{|A|,|V(X_i)\setminus A|\}<|\partial A|$ for any $i$
  and $A\subseteq V(X_i)$.} \]

  Margulis gave the first explicit construction of expanders. He made
  the fundamental observation that the Cayley graphs of finite quotients
  of a discrete group with Kazhdan's property~(T) form a family of
  expanders. Based on his ideas, Selberg's 3/16-theorem implies
  that $\{{\rm Cay}(\pi_m({\rm SL}_2(\ZZZ)),\pi_m(\Omega))\}_m$ is a
  family of expanders, where
  \[
    \Omega:=
    \left\{\left[\begin{array}{cc}
      1&\pm 1\\
      0&1
      \end{array}\right],
    \left[\begin{array}{cc}
      1&0\\
      \pm1&1
    \end{array}\right]\right\}
  , \]
  $m$ runs through all the positive integers, and $\pi_m$ is the
  reduction map modulo $m$. Many others (Burger, Sarnak and Clozel, to name a
  few) have studied the analytic
  behavior of congruence quotients of arithmetic lattices  using automorphic forms and
  representation theory. Their work resulted in the following.

  \begin{theorem}\label{t:PropertyTauArithmetic}
    Let $\mathbb{G}\subseteq \mathbb{GL}_n$ be a semisimple
    $\QQQ$-group and $\Gamma:=\mathbb{G}\cap {\rm GL}_n(\ZZZ_S)$,
    where $S$ is a finite set of primes. Assume $\Gamma$ is an infinite
    group that is generated by a finite (symmetric) set
    $\Omega=\Omega^{-1}$. Then the Cayley
    graphs ${\rm Cay}(\pi_m(\Gamma),\pi_m(\Omega))$ form a family
    of expanders as $m$ runs through positive integers.
  \end{theorem}

  Lubotzky was the first to ask if Theorem~\ref{t:PropertyTauArithmetic}
  holds for a {\em thin group}. Specifically he asked
  if $\{{\rm Cay}(\pi_p(\Gamma),\pi_p(\Omega))\}_p$ is a family of
  expanders, where
  \(
    \Gamma=
      \langle\Omega\rangle
  \),
  \[
    \Omega=\left\{\left[
      \begin{array}{cc}
        1&\pm 3\\0&1
      \end{array}\right],
    \left[
      \begin{array}{cc}
        1&0\\\pm3&1
      \end{array}
    \right]\right\}
  , \]
  and $p$ runs through all the primes. This question essentially asks
  if the {\em analytic behavior} of the congruence quotients of
  a linear group (under suitable conditions) is dictated by its
  Zariski topology.

  In a groundbreaking work based on a result of Helfgott~\cite{Gol-Hel},
  Bourgain and Gamburd~\cite{Gol-BG} answered Lubotzky's question
  affirmatively. These works were the starting point of a chain of
  fundamental results on the expansion properties of linear
  groups and their applications in other branches of mathematics
  (see the surveys~\cites{Gol-Lub, Gol-Kow, Gol-MSRI}).
    The following result was an essential part of the proof of the fundamental
  theorem of affine sieve~\cite{Gol-SGS}.

  \begin{theorem}[Golsefidy-Varj\'u \cite{Gol-SGV}]
    Let $\Omega\subseteq {\rm GL}_n(\ZZZ_S)$ be a finite symmetric set
    and $\Gamma=\langle\Omega\rangle$.
    Then $\{{\rm Cay}(\pi_q(\Gamma),\pi_q(\Omega))\}_q$ is a family of
    expanders as $q$ runs through square-free integers if and only if
    the Zariski connected component $\mathbb{G}^{\circ}$ of the
    Zariski closure $\mathbb{G}$ of $\Gamma$ is perfect, i.e.
    $\mathbb{G}^{\circ}=[\mathbb{G}^{\circ},\mathbb{G}^{\circ}]$.
  \end{theorem}

\endgroup

\begingroup
  \newcommand{\T}{\mathcal T}
  \newcommand{\Mod}{\mathrm {Mod}}
  \newcommand{\Diff}{\mathrm {Diff}}
  \newcommand{\Out}{\mathop{\mathrm{Out}}}
\recordedtalk{Outer automorphisms of free groups and tropical geometry}
     {\pers{Lizhen~Ji}{University of Michigan}}
     {Ji: Outer automorphisms of free groups and tropical geometry}
  Let $F_n$ be the free group on $n$ generators (with $n\geq 2$),
  and let $\Out(F_n)=\mathrm{Aut}(F_n)/\mathrm{Inn}(F_n)$ be its
  outer automorphism group. The group $\Out(F_n)$ is one of
  the most basic groups in combinatorial group theory, and it has been
  extensively studied.

  One key tool for understanding the properties
  of $\Out(F_n)$ is its action on a certain space~$X_n$, called \emph{outer
  space}. By definition (see \cite{CullerVogtmann}), $X_n$~is the space of equivalence classes of marked metric graphs
  whose fundamental group is~$F_n$.
  (A \emph{metric} on a graph~$\Gamma$ is an assignment of lengths to the edges of the graph, such that
  the sum is~$1$. A \emph{marking} of~$\Gamma$ is a homotopy equivalence from~$\Gamma$ to the wedge of
  $n$~circles. Roughly speaking, this is an identification of $\pi_1(\Gamma)$ with~$F_n$, but no basepoint has been fixed,
  so the identification is only well-defined up to an inner automorphism.
  Equivalence classes are taken with respect to a natural notion of isomorphism of marked metric graphs.)
  Since $\Out(F_n)$ acts on~$X_n$ (by changing the marking),
  information about~$X_n$ can yield information about $\Out(F_n)$.

  According to the celebrated Erlangen program of Klein,
  an essential part of the geometry of a space is concerned with
  invariants of the isometries (or symmetries) of the space.
  Similarly, an essential part of the geometry of a group is
  to find and understand spaces on which the group acts and
  preserves suitable additional structures of the space.
  It is often natural to require the space to be a metric
  space and the action to be by isometries.

  In classical geometry,  a metric space is a complete Riemannian
  manifold such as the Euclidean space, the sphere, or the
  hyperbolic space. But it is also important to consider metric
  spaces that are not manifolds. Examples include Tits buildings
  and Bruhat-Tits buildings for real and $p$-adic semisimple
  Lie groups. These are simplicial complexes with a natural
  complete Tits metric so that the groups act isometrically
  and simplicially on them. Of course, the rich combinatorial
  structure also makes their geometry interesting.

  The outer space $X_n$  can be realized in a natural way as a subset of a finite-dimensional
  simplicial complex (with infinitely many simplices) and hence admits a natural simplicial
  metric $d_0$. This metric is invariant under $\Out(F_n)$, 
  but it is not complete, because $X_n$ is missing the faces of some simplices.
  Since the completeness of metrics is a basic condition
  that is important for many applications, the following problem is very natural 
   (see  \cite[Question 2]{Ji-briv} for more
  discussion).

  \begin{problem}\label{prob-1.1}
  Construct  complete geodesic metrics on~$X_n$ that are invariant under $\Out(F_n)$.
  \end{problem}

  As explained above, solutions to this problem provide
  geometries of $\Out(F_n)$ in a certain sense.
  Another motivation for this problem comes from the analogy
  with other important groups in geometric group theory:
  arithmetic subgroups $\Gamma$ of semisimple Lie groups
  $G$ and mapping class groups $\Mod_{g,n}$ of surfaces
  of genus $g$ with $n$ punctures. \emph{A lot of work
  on $\Out(F_n)$ is motivated by results obtained
  for these families of groups.}

  An arithmetic group $\Gamma$ acts on the symmetric space $G/K$,
  and the mapping class group $\Mod_{g,n}$ acts on the
  Teichm\"uller space $\T_{g,n}$ of Riemann surfaces of genus
  $g$ with $n$ punctures. The symmetric space $X=G/K$ admits a
  complete $G$-invariant Riemannian metric and $\Gamma$ acts
  isometrically and properly on $X$. The Teichm\"uller space
  $\T_{g,n}$ admits several complete Riemannian and Finsler
  metrics, such as the Bergman and Teichm\"uller metrics,
  and $\Mod_{g,n}$ acts isometrically and properly on them.

  Though $X_n$ is not a manifold, it is a locally finite
  simplicial complex and hence is a canonically stratified space
  with a smooth structure. Therefore, the following problem
  also seems to be  natural in view of the above analogy.

  \begin{problem}\label{prob-1.2}
    Construct  piecewise smooth Riemannian metrics on~$X_n$ that are invariant under
    $\Out(F_n)$ and whose induced length metrics
    are complete geodesic metrics. Furthermore,
    the quotient $\Out(F_n)\backslash X_n$ has
    finite volume.
  \end{problem}

  With respect to the smooth structure on $X_n$ as a
  canonically stratified space, it is natural to require
   the Riemannian metric on $X_n$ to be smooth in the
  sense of stratified spaces. 
  Problems~\ref{prob-1.1} and~\ref{prob-1.2} are solved by
  the following theorem.

  \begin{theorem}[Ji \cite{Ji-ji}]
    There exist several explicitly constructed complete
    geodesic metrics and complete piecewise-smooth Riemannian
    metrics on $X_n$  that are invariant under $\Out(F_n)$.
  \end{theorem}

The proof of this theorem utilizes tropical geometry, which is 
algebraic geometry over the tropical
  semifield. (This is a rapidly developing subject --- see
  \cites{Ji-mi,Ji-ims} and the references therein.)
  The theory is applicable because metric graphs
  can be identified with tropical curves.
  There is a tropical Jacobian
  map from the moduli space of tropical curves to the moduli
  space of  principally polarized tropical abelian varieties, and
  the desired metrics are obtained by combining this map 
  with  the simplicial metric~$d_0$ on~$X_n$. This
  application of tropical geometry to
  geometric group theory might be of independent interest.
  
  Now that complete invariant geodesic metrics have been constructed
  on~$X_n$, one natural problem is to understand
  how these metrics can be used to study $X_n$ and $\Out(F_n)$.
 Their construction might be the first step
  towards a metric theory of outer space
  \cite[Question 1.2]{Ji-briv}.
\endgroup

\begingroup
  \newcommand{\RootSystem}{\Phi}
  \newcommand{\LieGroup}{G}
  \newcommand{\TheLattice}{\Gamma}
  \newcommand{\IrrReps}{r}
  \newcommand{\NumReps}{R}
  \newcommand{\GroupScheme}{\mathbf{G}}
  \newcommand{\PrimeSet}{S}
  \newcommand{\OkaRing}{\mathcal{O}}
  \newcommand{\SInts}{\OkaRing_{\PrimeSet}}
  \newcommand{\NumField}{k}
  \newcommand{\ThePlace}{v}
  \newcommand{\PlaceRing}{\OkaRing_{\ThePlace}}
  \newcommand{\SL}{\operatorname{SL}}
  \newcommand{\GL}{\operatorname{GL}}
  \newcommand{\SU}{\operatorname{SU}}
  \newcommand{\rmod}{/}
\recordedtalk{Representation growth of arithmetic groups}
     {\pers{Benjamin~Klopsch}{Heinrich Heine University D{\"u}sseldorf}
      and \\
      \pers{Christopher~Voll}{Bielefeld University}}
      {Klopsch and Voll: Representation growth of arithmetic groups}
  The talks of B.\,Klopsch and C.\,Voll were coordinated to essentially
  constitute another mini-course, so they are summarized here as a single unit.
  
  For a group $\TheLattice$, let $\IrrReps_n(\TheLattice)$ denote the number
  of irreducible complex representations (up to isomorphism) in dimension~$n$. 
  The group $\TheLattice$ is called \notion{representation rigid} if
  $\IrrReps_n(\TheLattice)$ is finite for every~$n$.

  A rigid group $\TheLattice$ has
  \notion{polynomial representation growth} {\small (PRG)}
  if $\IrrReps_n(\TheLattice)$ grows at most polynomially. Equivalently,
  one can ask that the partial sums
  \(
    \NumReps_n(\TheLattice)
    :=
    \sum_{i=1}^n \IrrReps_i(\TheLattice)
  \)
  grow at most polynomially. In this case, it is profitable to
  encode the numbers $\IrrReps_n(\TheLattice)$ into the
  \notion{(representation) zeta function} of $\TheLattice$:
  \[
    \zeta_{\TheLattice}(s)
    :=
    \sum_{n=1}^\infty \frac{ \IrrReps_n(\TheLattice) }{ n^s }
  .\]
  The {\small PRG}-condition ensures that $\zeta_{\TheLattice}$
  converges absolutely on some complex right half-plane. Conversely, the abscissa
  of convergence determines the rate of growth of the sequence
  $\{\NumReps_n(\TheLattice)\}$. Computing this abscissa is therefore a central problem in the subject
  of representation growth of groups. In particular, it is important to understand how the abscissa varies (and how
  other analytic invariants of~$\zeta_{\TheLattice}$ vary), as $\TheLattice$ ranges over interesting
  classes of groups.

  Consider the case that $\TheLattice$ is an arithmetic group in
  characteristic zero. For simplicity, assume that
  $\TheLattice=\GroupScheme(\SInts)$, where $\GroupScheme$ is a connected,
  simply connected, semisimple algebraic group defined over a number field $\NumField$
  with ring of $\PrimeSet$-integers $\SInts$.
  \begin{theorem}[Lubotzky-Martin \cite{MR2121543}]
    The group
    $\TheLattice$ has {\small PRG} if and only if it has the
    Congruence Subgroup Property {\small (CSP)}.
  \end{theorem}
  From now on we assume, again for simplicity, that the congruence kernel
  of $\TheLattice$ is trivial.
  In this case, the zeta function has an Euler product decomposition
  that is a consequence of Margulis superrigidity:
  \begin{proposition}[Larsen-Lubotzky \cite{MR2390327}]
    If $\TheLattice$ has trivial congruence kernel then the zeta
    function $\zeta_\TheLattice$ has an Euler product decomposition:
    \begin{equation}\label{eq:euler}
      \zeta_{\TheLattice}(s)
      =
      \zeta_{\GroupScheme(\CCC)}(s)^{[\NumField:\QQQ]}
      \prod_{\ThePlace\not\in\PrimeSet}{
        \zeta_{\GroupScheme(\PlaceRing)}
      }
    \end{equation}
  \end{proposition}
  The archimedean factor
  \(
    \zeta_{\GroupScheme(\CCC)}
  \)
  enumerates \emph{rational} representations of the algebraic group
  $\GroupScheme(\CCC)$. This factor, known as the Witten zeta function,
  is comparatively well understood in terms of highest weight theory and
  the Weyl character formula. The non-archimedean factors
  \(
    \zeta_{\GroupScheme(\PlaceRing)}
  \)
  enumerate \emph{continuous} representations of the $p$-adic analytic
  groups $\GroupScheme(\PlaceRing)$,
  and there exists a well-developed Lie theory for their
  principal congruence subgroups $\GroupScheme^m(\PlaceRing)$.
  In particular, the Kirillov orbit method sets up a
  $1$-$1$-correspondence between continuous, irreducible complex
  representations of these pro-$p$~groups and finite co-adjoint
  orbits in the dual of their Lie algebras. Deligne-Lusztig
  theory governs the representation theory of the finite groups
  of Lie type
  $\GroupScheme(\PlaceRing) \rmod \GroupScheme^1(\PlaceRing)$.

  Avni, Klopsch, Onn, and Voll have recently created a framework for the
  study of the local representation zeta functions $\zeta_{\GroupScheme(\PlaceRing)}$,
  with a view toward analyzing their Euler products.
 Formulas for the zeta functions of pro-$p$-groups of the form
  $\GroupScheme^m(\PlaceRing)$ can be obtained by developing methods from
  $p$-adic integration. These formulas yield, in particular, proofs of
  \emph{local functional equations} upon inversion of the residue field
  characteristic \cite{KuV-AKOV1}.

  For analysis of the analytic properties of Euler products such as
  \eqref{eq:euler}, control over the zeta functions of congruence subgroups
  is not sufficient. In \cite{KuV-AKOV4}, powerful machinery from model theory
 (viz., integrals of quantifier-free definable functions) is developed to
  ``approximate'', uniformly over large sets of places, the Clifford
  theory connecting the representation theory of the groups
  $\GroupScheme(\PlaceRing)$ with the representation
  theory of their congruence subgroups. This approach requires new
  insights into the behavior of representation growth of groups such as
  $\GroupScheme(\OkaRing)$ under base change, combined with a detailed
  analysis of representation zeta functions of finite groups of Lie type.
  This yields, in particular, a proof of the following theorem, which implies that the degree of
  representation growth is invariant under ring extensions.
  \begin{theorem}[Avni-Klopsch-Onn-Voll \cite{KuV-AKOV4}]
    For every irreducible root system $\RootSystem$, there is a constant
    $\alpha_\RootSystem$ that equals the abscissa of convergence of the
    representation zeta function for any group $\TheLattice=\GroupScheme(\SInts)$
    satisfying the {\small CSP}, such that
    $\GroupScheme$ has absolute root system $\RootSystem$.
  \end{theorem}
  This is related to the following conjecture, which is a refinement of Serre's
  conjecture on the Congruence Subgroup Property of lattices in higher-rank
  groups:
  \begin{conjecture}[Larsen-Lubotzky \cite{MR2390327}]
    Let $\LieGroup$ be a higher-rank semisimple locally compact group and let
    $\TheLattice_1$ and $\TheLattice_2$ be two irreducible lattices
    in $G$. Then the corresponding representation zeta functions have
    the same abscissa of convergence.
  \end{conjecture}
  For groups that have the {\small CSP} (which is required by the theorem, but not by the conjecture), the theorem's hypothesis is weaker  than the conjecture's hypothesis.
  Namely, the conjecture requires $\Gamma_1$ and~$\Gamma_2$ to be contained in a common ambient group~$G$, but the theorem
  only requires the ambient groups to have the same absolute root system.
  Further details can be found in the survey \cite{KuV-survey}.
\endgroup

\begingroup
\talk{Representation and characteristically equivalent arithmetic lattices}
     {\pers{C.~S.~Rajan}{Tata Institute of Fundamental Research}}
     {Rajan: Representation and characteristically equivalent arithmetic lattices}
  The inverse spectral problem is to recover the properties of a compact
  Riemannian manifold~$M$ from the knowledge of the spectrum of the Laplace operator
  (or of a more general Laplacian type operator) acting on the space of smooth
  functions  on~$M$.
  It is known, for example,
  that the spectra on functions determines the dimension, volume and
  the scalar curvature of $M$.

  Examples of non-isometric compact Riemannian manifolds
  which are isospectral on functions have been given by Milnor in the
  context of flat tori, and by  Vigneras for compact hyperbolic surfaces
  \cite{Raj-V}. Sunada gave a general method in analogy
  with a construction in arithmetic  \cite{Raj-S}.

  In many of these constructions, the manifolds are quotients by
  finite groups of a fixed Riemannian manifold. The question arises
  whether isospectral manifolds are indeed commensurable, i.e., have a
  common finite cover. In the context of
  Riemannian locally symmetric spaces this question has been studied
  by various authors \cites{Raj-R,Raj-CHLR,Raj-PR,Raj-LSV}
  assuming that the spaces are isospectral  for the
  Laplace-Beltrami operator acting on functions.
  In \cite{Raj-PR}, Gopal Prasad and A. S. Rapinchuk address this question
  in full generality, and get conditional commensurability type results for
  isospectral, compact locally symmetric spaces. For this when the
  locally symmetric spaces are of rank at least two, they have to
  assume the validity of Schanuel's conjecture on transcendental
  numbers. Another hypothesis they are required to make is that the base
  field is totally real and the group is anisotropic at all but one real
  place.

  Chandrasheel~Bhagwat, 
  Supriya~Pisolkar, and C.~S.~Rajan \cite{Raj-BPR}
  considered this question,
  assuming the stronger hypothesis
  that the lattices defining the locally symmetric spaces
  are representation equivalent, rather than isospectral on
  functions (see \cite{Raj-DG}).
  They were able to obtain unconditionally similar conclusions as in
  \cite{Raj-PR}
  for representation equivalent lattices, for example
  without invoking Schanuel's conjecture, and
  also extend the application to representation equivalent
  $S$-arithmetic lattices. In the process, they introduced a new relation of
  characteristic equivalence of lattices, stronger than weak
  commensurability.

  The proofs are distilled from the arguments given in \cite{Raj-PR}. The
   stronger hypothesis  simplifies  the  arguments used in
  \cite{Raj-PR}.
\endgroup

\begingroup
\recordedtalk{On the congruence subgroup problem}
     {\pers{Andrei~Rapinchuk}{University of Virginia}}
     {Rapinchuk: On the congruence subgroup problem}
  The talk was a brief survey of, and a progress report on, the
  congruence subgroup for algebraic groups over global fields. Let $G$
  be an absolutely almost simple algebraic group defined over a global
  field $K$, and let $S$ be a (not necessarily finite) set of places
  containing all the archimedean ones when $K$ is a number field. Then one
  considers the completions $\widehat{G}^S$ and $\overline{G}^S$ of
  the group $G(K)$ of rational point with respect to the
  $S$-arithmetic and the $S$-congruence topologies (see
  \cites{AR-PR-Milnor,AR-Rag-CSP} for precise definitions), and
  defines the $S$-{\it congruence kernel} $C^S(G)$ to be the kernel of
  the natural continuous surjective homomorphism $\widehat{G}^S
  \stackrel{\pi}{\longrightarrow} \overline{G}^S$. The congruence
  subgroup problem in this situation is the question about the
  computation of $C=C^S(G)$. The main conjecture, due to Serre, states
  that {\it $C$ should be \emph{finite} if $\mathrm{rk}_S\: G :=
  \sum_{v \in S} \mathrm{rk}_{K_v}\: G$ is $\geqslant 2$ and
  $\mathrm{rk}_{K_v}\: G > 0$ for all nonarchimedean $v \in S$, and
  \emph{infinite} if $\mathrm{rk}_S\: G = 1$} (here
  $\mathrm{rk}_{K_v}\: G$ denotes the rank of $G$ over the completion
  $K_v$). The talk focused primarily on the higher rank case of
  Serre's conjecture (the structure of $C$ in the rank one situation
  has been determined in many cases by O.~V.~Mel'nikov, A.~Lubotzky and
  P.~A.~Zalesskii with various co-authors - see the references in
  \cite{AR-PR-Milnor}). First, it was explained that modulo the
  Margulis-Platonov conjecture on the structure of normal subgroups of
  $G(K)$, which has already been proved in the majority of cases, the
  finiteness of $C$ is equivalent to its {\it centrality}, i.e., to the
  fact that it lies in the center of $\widehat{G}^S$, in which case
  $C$ (or more precisely, its Pontrjagin dual) is isomorphic to the
  {\it metaplectic kernel} $M(S , G)$ (see \cite{AR-PR-Milnor},
  \cite{AR-Rag-CSP} for the details). Second, the effort to compute $M(S
  , G)$ initiated by C.~Moore and continued by Matsumoto, Deodhar,
  Prasad-Raghunathan and others was completed in \cite{AR-PR-Met}. In
  essence, the final result says that $M(S , G)$ is always finite and
  is isomorphic to a subgroup of the group $\mu_K$ of roots of unity
  in $K$ (in particular, it is always a finite cyclic group), and in
  fact is trivial under some rather general additional assumptions.
  So, the focus in the higher rank case of the congruence subgroup
  problem is currently on finding a general approach to the proof of
  centrality. The centrality has been established in a number of cases
  using a variety of techniques (cf.\ \cite{AR-PR-Milnor,AR-Rag-CSP}), but there are still anisotropic groups (e.g.,
  $K$-groups of the form $\mathrm{SL}_{1 , D}$ where $D$ is a
  finite-dimensional central division algebra over $K$) that defy all
  efforts. The talk discussed some approaches to proving the
  centrality that do not require any case-by-case considerations. One
  of the approaches relies on the analysis of the centralizers of
  elements from $G(K)$ in $\widehat{G}^S$. To formulate a result in
  this direction, recall that, without loss of generality, we may
  assume that $\mathrm{rk}_S\: G > 0$, and then the $S$-congruence
  completion $\overline{G}^S$ can be identified with the group of
  $S$-adeles $G(\mathbb{A}(S))$ by the strong approximation theorem
  (see \cite{AR-R-SA} for a recent survey on strong approximation).

  \begin{proposition}\label{AR-prop-1}
    Assume that there exists $n
    \geqslant 1$ such that for any regular semi-simple element $t \in
    G(K)$ and its centralizer $T = Z_G(t)$ we have
    $$
    \pi(Z_{\widehat{G}^S}(t)) \supset T(\mathbb{A}(S))^n.
    $$
    Then $C$ is central.
  \end{proposition}

  If we let $\widehat{\ }$ and $\bar{\parbox[t]{3mm}{\ } }$ denote the
  closure in $\widehat{G}^S$ and $\overline{G}^S$, respectively, then
  clearly, $Z_{\widehat{G}^S}(t)$ contains $\widehat{T(K)}$, and since
  $\pi(\widehat{T(K)}) = \overline{T(K)}$, we obtain the following.

  \begin{corollary}\label{AR-cor-1}
    If there exists $n \geqslant 1$
    such that
    $$
    \overline{T(K)} \supset T(\mathbb{A}(S))^n
    $$
    for any $K$-torus $T$ in $G$ then $C$ is central.
  \end{corollary}

  Thus, the centrality would be a consequence of the property of
  almost strong approximation in all $K$-tori of $G$ for a given $S$.
  The bad news is that this property \emph{never} holds for a
  nontrivial $K$-torus if $S$ is finite. More precisely, one uses the
  fact that $T$ admits coverings $T \to T$ of any degree to show that
  the quotient $T(\mathbb{A}(S))/\overline{T(K)}$ has infinite
  exponent (cf. \cite{AR-R-SA}). But there is a little bit of good news
  at the other end of the spectrum, viz. when $S$ is co-finite, i.e.
  $S = V^K \setminus S_0$ where $S_0$ is a finite set of
  nonarchimedean places. In this case,
  $$
  T(\mathbb{A}(S)) = T_{S_0} := \prod_{v \in S_0} T(K_v),
  $$
  and one shows that $T_{S_0}/\overline{T(K)}$ has finite exponent
  which can be bounded by a function depending only on $\dim T$.

  \begin{corollary}[Semi-local case]\label{AR-cor-2}
    If $S$ is co-finite then $C$ is central (in fact, trivial).
  \end{corollary}

  In fact, the property of almost strong approximation holds in tori
  if $S$ almost contains a generalized arithmetic progression (subject
  to some natural assumptions), and one can formulate the
  corresponding result for the centrality of $C$.

  \vskip2mm

  Corollary~\ref{AR-cor-2} leads to another condition for the
  centrality of $C$. To
  formulate it, we observe that using the identification
  $\overline{G}^S = G(\mathbb{A}(S))$, one can think of $G(K_v)$ for
  any $v \notin S$ as a subgroup of $\overline{G}^S$.

  \begin{theorem}\label{AR-thm-1}
    Assume that for each $v \notin S$
    there exists a subgroup $\mathcal{G}_v \subset \widehat{G}^S$ so
    that
    \begin{enumerate}
      \item $\pi(\mathcal{G}_v) = G(K_v)$;
      \item $\mathcal{G}_{v_1}$ and $\mathcal{G}_{v_2}$ commute
            elementwise for any $v_1 \neq v_2$;
      \item the subgroup generated by the $\mathcal{G}_v$'s is dense
            in $\widehat{G}^S$.
    \end{enumerate}
    Then $C$ is central.
  \end{theorem}

  This result can be used to give a relatively short proof of Serre's
  conjecture for $K$-isotropic groups and also for $K$-anisotropic
  groups of exceptional types $\textsf{E}_7, \textsf{E}_8$ and
  $\textsf{F}_4$ (recall that these types split over a quadratic
  extension of $K$).

  Finally, the ideas involved in the proof of Theorem~\ref{AR-thm-1}
  can be used to
  provide some information about $C$ in the rank one case. First,
  recall the following general dichotomy (assuming the truth of the
  Margulis-Platonov conjecture): {\it $C$ is either central and
  finite, or is \emph{not} finitely generated} (e.g., the congruence
  kernel for the group $\mathrm{SL}_2(\ZZZ)$ is known to be the
  free profinite group of countable rank). Nevertheless, Lubotzky
  proved that if $\Gamma = G(\mathcal{O}(S))$ is the corresponding
  $S$-arithmetic subgroup then $C$ is {\it finitely generated} as a
  normal subgroup of $\widehat{\Gamma}$ (this is a consequence of
  finite presentation of the group of integral ideles
  $\overline{\Gamma}$). G.~Prasad and A.~Rapinchuk showed that in some
  situations $C$ is virtually generated by a {\it single} element as a
  normal subgroup of $\widehat{G}^S$. (Note that, between this
  result and that of Lubotzky, neither one implies the other.)

  \begin{theorem}\label{AR-thm-2}
    Assume that $K$ is a number field
    and $G$ is $K$-isotropic. Then there exists $c \in C$ such that if
    $D \subset C$ is the closed normal subgroup of $\widehat{G}^S$
    generated by $c$ then $C/D$ is a finite cyclic group.
  \end{theorem}
  (For $G = \mathrm{SL}_2$ this element $c$ can be written down
  explicitly.)
\endgroup

\begingroup
\recordedtalk{On the conjecture of Borel and Tits for abstract
homomorphisms of algebraic groups}
     {\pers{Igor~Rapinchuk}{Yale University}}
     {Rapinchuk: On the conjecture of Borel and Tits for abstract
homomorphisms of algebraic groups}
  The general philosophy in the study of abstract homomorphisms
  between groups of rational points of algebraic groups is as
  follows. Suppose $\mathbf{G}$ and $\mathbf{G}'$ are algebraic groups that are
  defined over infinite fields $K$ and $K',$ respectively. Let
  \[
    \varphi \colon \mathbf{G}(K) \to \mathbf{G}'(K')
  \]
  be an abstract homomorphism between their groups of rational
  points. Then, under appropriate assumptions, one expects to be
  able to write $\varphi$ essentially as a composition
  $\varphi = \beta \circ \alpha,$ where
  $\alpha \colon \mathbf{G}(K) \to {}_{K'}\!\mathbf{G}(K')$ is induced by a
  \emph{field homomorphism}
  $\tilde{\alpha} \colon K \to K'$ (and $_{K'}\!\mathbf{G}$ is the group
  obtained from $\mathbf{G}$ by base change via $\tilde{\alpha}$), and
  $\beta \colon _{K'}\!\mathbf{G}(K') \to \mathbf{G}'(K')$ arises from a
  $K'$-defined \emph{morphism of algebraic groups}
  $_{K'}\!\mathbf{G} \to \mathbf{G}'$. Whenever $\varphi$ admits such a
  decomposition, one generally says that it has a
  \notion{standard description}.

  The following conjecture
  of Borel and Tits  is a major open question. Recall that, for an
  algebraic group~$\mathbf{G}$ defined over a field $k$, one denotes by $\mathbf{G}^+$
  the subgroup of $\mathbf{G}(k)$ generated by the $k$-points of split
  (smooth) connected unipotent $k$-subgroups.

  \begin{conjecture}[Borel-Tits {\cite[8.19]{Rap-BT}}]\label{Rap-conj-bt}
    Let $\mathbf{G}$ and $\mathbf{G}'$ be algebraic groups defined over infinite
    fields $k$ and $k'$, respectively.
    If $\rho \colon \mathbf{G}(k) \to \mathbf{G}'(k')$ is any abstract homomorphism,
    such that $\rho(\mathbf{G}^+)$ is Zariski-dense in $\mathbf{G}'(k')$, then
    there exists a commutative finite-dimensional
    $k'$-algebra~$B$ and a ring homomorphism
    $f \colon k \to B$, such that
    \[
      \rho \vert_{\mathbf{G}^+} = \sigma \circ r_{B/k'} \circ F
    , \]
    where 
    	\begin{itemize}
	\item $F \colon \mathbf{G}(k) \to {}_{B}\mathbf{G}(B)$ is induced by
    $f$ (and $_{B}\mathbf{G}$ is obtained by change of
    scalars),
    \item
    $r_{B/k'} \colon {}_{B}\mathbf{G}(B) \to R_{B/k'} (_{B}\mathbf{G})(k')$
    is the canonical isomorphism  (here $R_{B/k'}$ denotes the
    functor of restriction of scalars), and 
    \item $\sigma$ is a
    rational $k'$-morphism of $R_{B/k'} (_{B}\mathbf{G})$ to $\mathbf{G}'$.
    \end{itemize}
  \end{conjecture}

  In their fundamental paper \cite{Rap-BT}, Borel and Tits proved
  the conjecture for $\mathbf{G}$  an absolutely almost simple $k$-isotropic
  group and $\mathbf{G}'$ a reductive group. Shortly after the conjecture
  was formulated, Tits \cite{Rap-T} sketched a proof of~\ref{Rap-conj-bt}
  in the case that $k = k' = \RRR.$ Prior to the recent work of I.~Rapinchuk that is described below, the
  only other available result was due to L.~Lifschitz and
  A.~S.~Rapinchuk \cite{Rap-LR}, where the conjecture was essentially
  proved in the case where $k$ and $k'$ are fields of
  characteristic~$0$, $\mathbf{G}$ is a universal Chevalley group, and
  $\mathbf{G}'$ is an algebraic group with commutative unipotent radical.

  While the above results only deal with abstract homomorphisms
  of groups of points over \emph{fields}, it should be pointed
  out that there has also been considerable interest and activity
  in analyzing abstract homomorphisms of higher rank arithmetic
  groups and lattices (e.g., the work of Bass, Milnor, and
  Serre \cite{Reid-BMS} on the congruence subgroup problem and
  Margulis's Superrigidity Theorem \cite[Chap. VII]{Rap-Mar}).
  However, relatively little was previously known about abstract
  homomorphisms of groups of points over \emph{general commutative
  rings}, which has been the primary focus of I.~Rapinchuk's work in this area.

  To state the new results, we first need to fix some notation.
  Let $\Phi$ be a reduced irreducible root system of rank $\geq 2$
  and $\mathbf{G}$ be the corresponding universal Chevalley-Demazure group
  scheme over $\ZZZ.$ For any commutative ring $R$, we
  denote by $\mathbf{G}(R)^+$ the subgroup of $\mathbf{G}(R)$ generated by the
  $R$-points of the canonical one-parameter root subgroups (usually
  called the {\it elementary subgroup}). 
  
  The first theorem below is a
  rigidity result for abstract representations
  \[
    \rho \colon \mathbf{G}(R)^+ \to GL_n (K),
  \]
  where $K$ is an algebraically closed field. In its statement, 
  for a finite-dimensional commutative $K$-algebra $B$,
  the group of rational points $\mathbf{G}(B)$ is viewed as an algebraic
  group over $K$ by using the functor of restriction of scalars.
  Furthermore, given a commutative ring $R$, we will say that
  $(\Phi, R)$ is a \emph{nice pair} if $2 \in R^{\times}$
  whenever $\Phi$ contains a subsystem of type $B_2$,
  and $\{2, 3\} \subseteq R^{\times}$ if $\Phi$ is of type $G_2.$

  \begin{theorem}[I.~Rapinchuk {\cite[Main Theorem]{Rap-IR}}]\label{rap:main}
    Let $\Phi$ be a reduced irreducible root system of rank
    $\geq 2$, $R$ a commutative ring such that $(\Phi, R)$ is
    a nice pair, and $K$ an algebraically closed field.
    Assume that $R$ is noetherian if
    $\mathrm{char} \: K > 0.$ Furthermore let $\mathbf{G}$ be the
    universal Chevalley-Demazure group scheme of type
    $\Phi$ and let $\rho \colon \mathbf{G}(R)^+ \to GL_n (K)$ be a
    finite-dimensional linear representation over $K$ of the
    elementary subgroup $\mathbf{G}(R)^+ \subset \mathbf{G}(R)$. Set
    $H = \overline{\rho (\mathbf{G}(R)^+)}$ (Zariski closure),
    and let $H^{\circ}$ denote the connected component of
    the identity of $H$. Then in each of the following
    situations
    \begin{enumerate}
      \item
        $H^{\circ}$ is reductive;
      \item
        $\mathrm{char} \: K = 0$ and $R$ is semilocal;
      \item
       $\mathrm{char} \: K = 0$ and the unipotent radical
       $U$ of $H^{\circ}$ is commutative,
    \end{enumerate}
    there exists a commutative finite-dimensional $K$-algebra
    $B$, a ring homomorphism $f \colon R \to B$ with
    Zariski-dense image, and a morphism
    $\sigma \colon \mathbf{G}(B) \to H$ of algebraic $K$-groups such that
    for a suitable subgroup $\Delta \subset \mathbf{G}(R)^+$ of
    finite index, we have
    \[
      \rho \vert_{\Delta} = (\sigma \circ F) \vert_{\Delta},
    \]
    where $F \colon \mathbf{G}(R)^+ \to \mathbf{G}(B)^+$ is the group homomorphism
    induced by $f$.
  \end{theorem}

  Thus, if $R = k$ is a field of characteristic $\neq 2$ or 3,
  then $R$ is automatically semilocal and $(\Phi, R)$ is a
  nice pair. Hence, Theorem 1 provides a proof of Conjecture~\ref{Rap-conj-bt}
  in the case that $\mathbf{G}$ is split and $K$ is an algebraically
  closed field of characteristic zero.

  Let us now describe some applications of Theorem~\ref{rap:main}
  to the study of character varieties of elementary subgroups of
  Chevalley groups. Let $K$ be an algebraically closed field of
  characteristic~0 and $R$ be a finitely generated commutative
  ring. As above, suppose that $\Phi$ is a reduced irreducible
  root system of rank $\geq 2$ and let $\mathbf{G}$ be the corresponding
  universal Chevalley-Demazure group scheme. Then the elementary
  subgroup $\mathbf{G}(R)^+$ has Kazhdan's property $(T)$ (see \cite{Rap-EJK}),
  hence is in particular a finitely generated group, and
  therefore, for any integer $n \geq 1$, one can consider the
  character variety $X_n (\Gamma)$.

  \begin{theorem}[I.~Rapinchuk \cite{Rap-IR1}]\label{rap:two}
    Let $\Phi$ be a reduced irreducible root system of rank
    $\geq 2$, $R$~a finitely generated commutative ring such
    that $(\Phi, R)$ is a nice pair, and $\mathbf{G}$ the universal
    Chevalley-Demazure group scheme of type $\Phi$.
    Denote by $\Gamma$ the elementary subgroup $\mathbf{G}(R)^+$ of
    $\mathbf{G}(R)$ and consider the $n$th character variety
    $X_n (\Gamma)$ of $\Gamma$ over an algebraically closed
    field $K$ of characteristic 0. Then there exists a
    constant $c = c(R)$ (depending only on $R$) such
    that $\kappa_n (\Gamma) := \dim X_n (\Gamma)$ satisfies
    \[
      \kappa_n (\Gamma) \leq c \cdot n
    \]
    for all $n \geq 1.$
  \end{theorem}

  The proof of Theorem~\ref{rap:two} makes extensive use of Theorem~\ref{rap:main}'s 
  description of the representations with non-reductive image.

  Another application of Theorem~\ref{rap:main} has to do with
  the problem of realizing complex affine varieties as character
  varieties of suitable finitely generated groups. This question
  was previously considered by M.~Kapovich and
  J.~Millson \cite{Rap-KM}, who showed that any affine variety $S$
  defined over $\QQQ$ is birationally isomorphic to an
  appropriate character variety of some Artin group $\Gamma.$
  By using Theorem~\ref{rap:main} with $K = \CCC$, it is possible to prove the
  following result.

  \begin{theorem}[I.~Rapinchuk \cite{Rap-IR2}]\label{rap:three}
    Let $S$ be an affine algebraic variety defined over
    $\QQQ.$ There exist a finitely generated group
    $\Gamma$ having Kazhdan's property $(T)$ and an integer
    $n \geq 1$ such that there is a biregular isomorphism
    of complex algebraic varieties
    \[
      S (\CCC) \to X_n (\Gamma) \setminus \{ [\rho_0] \},
    \]
    where $\rho_0$ is the trivial representation.
  \end{theorem}
\endgroup

\begingroup
  \newcommand{\SL}{\operatorname{SL}}
  \newcommand{\Sp}{\operatorname{Sp}}
  \newcommand{\PSL}{\operatorname{PSL}}
  \newcommand{\SO}{\operatorname{SO}}
\recordedtalk{All finite groups are involved in the mapping class group}
     {\pers{Alan~Reid}{University of Texas at Austin}}
    {Reid: All finite groups are involved in the mapping class group}
  Let $\Sigma_g$ be a closed orientable surface of genus $g\geq 1$, and
  $\Gamma\!_g$ its Mapping Class Group.

  A group $H$ is \emph{involved} in a group $G$ if there is a finite index subgroup
  $K<G$ so that $K$ subjects onto~$H$. The question as to whether
  every finite group is involved in a fixed~$\Gamma\!_g$ was raised by U.~Hamenst{\"a}dt in her talk at the 2009 Georgia Topology Conference.
  This is easily seen to hold for the case $g=1$ (since $\Gamma\!_1=\SL(2,{\bf
  Z})$ is virtually free) and for $g=2$ (since $\Gamma\!_2$ is large, see
 \cite{Reid-Ko}).
 In fact, it holds for all~$g$:

  \begin{theorem}[Masbaum-Reid \cite{Reid-MR}] \label{main}
    For all $g\geq 1$, every finite group is involved in $\Gamma\!_g$.
  \end{theorem}

  Although $\Gamma\!_g$ is well-known to be residually finite \cite{Reid-Gro},
  and therefore has a rich supply of finite quotients,
   very little seems known about what finite groups can arise as
  quotients of $\Gamma\!_g$ (or of subgroups of finite index), other than those finite quotients obtained from
  \[
    \Gamma\!_g \rightarrow \Sp(2g,\ZZZ) \rightarrow \Sp(2g,\ZZZ/N\ZZZ)
  .\]
  It should be emphasized that one
  cannot expect to prove Theorem \ref{main} by simply using
  the subgroup structure of the groups $\Sp(2g,\ZZZ/N\ZZZ)$. The reason
  for this is
  that, since $\Sp(2g,\ZZZ)$ has the Congruence Subgroup
  Property \cite{Reid-BMS}, it is well-known that not all finite groups
  are involved in
  $\Sp(2g,\ZZZ)$ (see \cite[Chapter 4.0]{Reid-LR}, for example).

  The main new idea in the proof of Theorem \ref{main} is to
  exploit the unitary representations arising
  in Topological Quantum Field Theory (TQFT), first constructed by
  Reshetikhin and
  Turaev \cite{Reid-RT}. (Actually, the proof uses the so-called $\SO(3)$-TQFT,
  following the skein-theoretical approach of \cite{Reid-BMVH} and the Integral
  TQFT refinement \cite{Reid-GM}.)

Since, as was mentioned above, the case $g = 1$ and the case $g = 2$ are easy, it suffices to deal with the case where $g\geq 3$.
Therefore, Theorem~\ref{main} easily follows from the next result,
  which gives
  many new finite simple groups of Lie type as quotients of~$\Gamma\!_g$.

  \begin{theorem}[Masbaum-Reid \cite{Reid-MR}] \label{main2}
    For each $g\geq 3$, there exist infinitely many $N$, such that,
    for each such $N$, there exist infinitely
    many primes~$q$, such that $\Gamma\!_g$ surjects onto the finite group $\PSL(N,\FFF_q)$,
    where $\FFF_q$ denotes the finite field of order~$q$.
  \end{theorem}

  In addition, \cite{Reid-MR} also shows that Theorem~\ref{main2}
  holds for the Torelli group (with $g\geq 2$).
\endgroup

\begingroup
\recordedtalk{Informal talk on Kac-Moody groups}
     {\pers{Bertrand~R{\'e}my}{Institut Camille Jordan}}
     {Remy: Informal talk on Kac-Moody groups}
%
  There are two kinds of Kac-Moody groups: the complete groups and the
  minimal groups. Both are discussed in J.~Tits' Bourbaki seminar
  \cite{Remy-Tits89}, and they both have the same
  algebraic origin, namely Kac-Moody Lie
  algebras \cite{Remy-Kac90}. (Kac-Moody 
  algebras are infinite-dimensional analogues of finite-dimensional reductive Lie algebras. In the classical finite-dimensional setting, the Serre presentation produces generators and relations for the Lie algebra from a Cartan matrix. In the infinite-dimensional setting, there is an analogous presentation that is produced from a generalized Cartan matrix.)
  The two kinds of groups also share a method of construction that imitates, in an
  infinite-dimensional context, the definition of Chevalley-Demazure
  group schemes for reductive groups \cite{Remy-Dem65}.
  Namely, the Kac-Moody group of either kind is a functor that is defined by a
  presentation. It is essentially a Steinberg presentation, generalizing an abstract presentation of
  the rational points of a split reductive group \cite{Remy-Tits87}.
  The main ingredients in the presentation are various
  completions of integral forms of (pieces of) universal enveloping
  algebras of Kac-Moody algebras \cite{Remy-Rem02}.

  Instead of
  working with the presentation of a Kac-Moody group, it is much more
  efficient to use its nice combinatorics (the existence of two twinned
  Tits systems).
  The geometric counterpart to these rich combinatorial properties is
  the existence of buildings on which the Kac-Moody group
  acts highly transitively.
  By definition, buildings are cell complexes that are the union of subcomplexes all
  isomorphic to a given Coxeter tiling; some additional incidence
  properties are required \cite{Remy-AB08}.
  They admit very useful metrics that are complete and non-positively
  curved. Moreover, non-spherical buildings are contractible, which
  suggests a fruitful analogy with symmetric spaces of non-compact type
  (an important tool in the study of Lie groups and their discrete
  subgroups).

  One of the valuable features of Kac-Moody theory is that it leads to
  intriguing new examples of groups.
  For example, the minimal Kac-Moody groups over finite fields provide infinitely
  many quasi-isometry classes of finitely presented simple groups \cite{Remy-CR10}.
  (The proof of simplicity has two main parts. First is the fact that non-central normal subgroups have finite index \cite{BaderShalom,Remy-induction}, which is analogous to a well-known result in the theory of arithmetic groups. Then the proof exploits a crucial ``weakly hyperbolic'' property of the geometry of non-Euclidean infinite Coxeter groups \cite{Remy-CR09}.)
  The maximal pro-$p$ subgroups in complete Kac-Moody groups are another class of groups that pose interesting challenges. Their first homology has just
  been computed \cite{Remy-CR13}, but their higher finiteness
  properties still need to be investigated.
\endgroup

\begingroup
\recordedtalk{Counting ends of rank one arithmetic orbifolds}
     {\pers{Matthew~Stover}{Temple University}}
     {Stover: Counting ends of rank one arithmetic orbifolds}
  If $N$ is a noncompact negatively curved locally symmetric
  space of finite volume, then it has a finite number of topological
  ends, or cusps. The following question, remarkably, remains wide
  open:
  \begin{question}
    Is there a one-cusped complete finite-volume hyperbolic
    $n$-manifold for every~$n$?
  \end{question}

  One-ended finite-volume orbifold quotients of hyperbolic $n$-space
  $\mathbf{H}^n$ are known for $n < 10$, but no one-cusped
  $n$-manifold is known for $n > 4$. For $n = 2, 3$, it is relatively
  easy to find examples with arithmetic fundamental group. For
  example, one can interpret the fact that the modular surface
  $\mathrm{PSL}_2(\ZZZ) \backslash \mathbf{H}^2$ has one cusp in
  terms of the fact that $\ZZZ$ is a principal ideal domain, and
  it is relatively easy to find a one-cusped manifold cover. The
  number of cusps of arithmetic hyperbolic 3-manifolds is closely
  related to the class number of an imaginary quadratic field k, which
  is an invariant that measures `how far' its ring of integers is from
  being a PID. 
  
  Given the ease with which one can build one-cusped
  arithmetic orbifolds in dimensions 2 and 3, one might hope that
  arithmetic techniques could provide one-cusped hyperbolic
  $n$-manifolds, or at least orbifolds, for all $n$. The following
  theorem shows that this is impossible.

  \begin{theorem}[Stover \cite{Stover}]\label{Sto-thm-1}
    One-cusped arithmetic hyperbolic $n$-orbifolds do not exist
    for any $n > 31$.
  \end{theorem}

  The proof relates the number of ends to the so-called class number
  of a certain quadratic form, then studies related number theoretic
  invariants that yield a lower bound on the number of ends. In
  fact, it gives an exact formula for the number of cusps for $\Gamma
  \backslash \mathbf{H}^n$ when $\Gamma$ is a natural generalization
  of the usual congruence subgroups of $\mathrm{PSL}_2(\ZZZ)$. The paper
  also constructs new one-cusped examples in dimensions 10 and 11.
 
  Theorem~\ref{Sto-thm-1} is actually a precise special case of the
  following much stronger finiteness theorem.

  \begin{theorem}[Stover \cite{Stover}]\label{Sto-thm-2}
    Fix $k > 0$. There are only finitely many commensurability
    classes of negatively curved arithmetic locally symmetric spaces
    that contain an element with $k$~ends.
  \end{theorem}

  The negatively curved locally symmetric spaces are hyperbolic
  $n$-space, complex hyperbolic $n$-space $\mathbf{H}_{\CCC}^n$,
  quaternionic $n$-space $\mathbf{H}_{\mathbb{H}}^n$, and the
  exceptional Cayley hyperbolic plane $\mathbf{H}_{\mathbb{O}}^2$. All
  finite-volume quaternionic hyperbolic $n$-orbifolds and Cayley
  hyperbolic $2$-orbifolds are arithmetic, so the arithmetic
  assumption in Theorem~\ref{Sto-thm-2} is superfluous and finiteness
  holds over all finite-volume quotients. In particular, for
  each $k > 0$, there is a
  constant $c_k$ such that finite-volume quaternionic hyperbolic
  $n$-orbifolds with $k$ ends do not exist for $n > c_k$.
\endgroup

\begingroup
  \newcommand{\Q}{\QQQ}
  \newcommand{\Z}{\ZZZ}
  \newcommand{\Sp}{\operatorname{Sp}}
\recordedtalk{Monodromy of arithmetic groups}
     {\pers{T.\ N.\,Venkataramana}{Tata Institute of Fundamental Research}}
     {Venkataramana: Monodromy of arithmetic groups}
  Let $f,g\in \Z[X]$ be
  polynomials of degree~$n$ that are monic with constant term one, and have
  no common root. Also assume that every root  of $fg$ is a root of unity, and  that $\{f,g\}$ is a ``primitive pair'' (see  \cite{Venk-SV}).
  
  \begin{theorem}[Beukers-Heckman \cite{Venk-BH}] 
  The  companion matrices  $A,B$  of $f,g$
  preserve  a nondegenerate  symplectic form~$\Omega $  on~$\Q  ^n$ and
  generate a Zariski dense  subgroup $\Gamma$ of the integral symplectic
  group  $\Sp_n(\Omega, \Z)$. 
  \end{theorem}

  It is known \cite{Venk-BH} that  $\Gamma$  is 
the monodromy group  of a suitable hypergeometric equation
  of type $_nF_{n-1}$.
  The following theorem determines  when the subgroup~$\Gamma$  has
  finite index in $\Sp_n(\Omega, \Z)$:  
  
  \begin{theorem}[Singh-Venkataramana \cite{Venk-SV}]
  If the  leading coefficient  of the polynomial  $f-g$ does  not exceed
  two, then $\Gamma$  is arithmetic. 
  \end{theorem}
  
  The method of  proof also shows that for
  the 14  examples of Calabi-Yau  threefolds listed in  \cite{Venk-CEYY}, the
  monodromy group is arithmetic.
\endgroup

\begingroup
  \newcommand{\GL}{\operatorname{GL}}
  \newcommand{\Z}{\ZZZ}
  \newcommand{\E}[1]{\mathrm{E}#1}
  \newcommand{\gd}{\operatorname{gd}}
  \newcommand{\pE}[1]{\underline{\mathrm{E}\!}\,#1}
  \newcommand{\pgd}{\underline{\smash{\operatorname{gd}}}}
  \newcommand{\pF}[1]{\underline{\mathrm{F}\!}_{\,#1}}
  \newcommand{\F}[1]{\mathrm{F}_{#1}}
\recordedtalk{Bredon finiteness properties of arithmetic groups}
     {\pers{Stefan~Witzel}{University of M{\"u}nster}}
     {Witzel: Bredon finiteness properties of arithmetic groups}
  Classifying spaces $\E{G}$ of groups $G$ have been studied for
  a long time. Its usefulness in
  understanding~$G$ partly depends on how finite $\E{G}$ is
  (or how finite it can be chosen to be). One natural measure of finiteness
  is dimension, which leads to the notion of the \emph{geometric
  dimension} $\gd G$ of the group~$G$. In a different direction,
  one investigates up to which dimension the action of $G$ on $\E{G}$
  is cocompact, which is encoded in the \emph{finiteness
  properties}~$\F{n}$.

  For groups with torsion, it is natural to allow the actions
  to have finite stabilizers. Thus, instead of classifying spaces $\E{G}$
  for free actions, one studies classifying spaces $\pE{G}$
  for proper actions. (All actions are by cell-permuting
  homeomorphisms on CW-complexes.) In the same way as for
  free actions, this gives rise to the notions of \emph{proper
  geometric dimension} $\pgd$ and \emph{proper finiteness
  properties} $\pF{n}$ (see \cite{MR0206946} and \cite{MR1027600}).

  For the study of classical (free) finiteness properties,
  there is a very useful criterion due to Brown \cite{MR885095}.
  If one can let $G$ act on a contractible space, in such a way
  that the stabilizers have good finiteness properties
  themselves, then the criterion relates the finiteness properties
  of $G$ to the essential connectivity of an orbit. (Essential
  connectivity is a technical property that measures how
  highly connected an orbit is in a coarse sense.) 
  
  Fluch and
  Witzel \cite{FluchWitzel2011} translated a homological version
  of Brown's criterion to actions with arbitrary families
  of stabilizers, in particular, to proper actions.
 Algebraically, proper
  finiteness properties correspond to finiteness properties
  of normalizers of finite groups. Topologically, they are
  reflected in the connectivity of fixed point sets of
  finite subgroups. The essential connectivity in this
  case is measured uniformly over all finite groups.

  A concrete family of examples illustrates how
  and why the classical finiteness properties of a group can differ
  from its proper finiteness
  properties. Namely, consider the stabilizers in $\GL_{n+1} \bigl( \Z[1/p] \bigr)$
  of two horospheres in the associated Bruhat--Tits building.
  In a special case, these groups were known  to be of type $\F{n-1}$,
  but not of type $\F{n}$, by work of Abels, Brown and others \cite{MR885096}.
  The classical finiteness properties generalize to the whole family.
  In contrast, the proper finiteness properties depend on the position of the
  horospheres. This is because the fixed-point sets of finite
  subgroups decompose as products of buildings, and the connectivity
  of horospheres in these depends on which direct factors are
  contained in one of the horospheres. In addition, the
  amount of torsion in the groups can also vary. A detailed analysis 
  yields the following theorem, which shows that the two types of finiteness
  properties can vary completely independently of each other.
  
  \begin{theorem}[Witzel \cite{Witzel2012}]
  For $0<m \le n$, there is a solvable algebraic group $\mathbf{G}$, such that, for every
odd prime~$p$, the group $\mathbf{G} \bigl (\ZZZ[1/p] \bigr)$ is
  	\begin{itemize}
	\item of type $\pF{m-1}$, but not $\pF{m}$, and, also, 
	\item of type $\F{n-1}$ but not $\F{n}$.
  \end{itemize}
  \end{theorem}
  \endgroup

\begingroup
\recordedtalk{Profinite topology on arithmetic groups}
     {\pers{Pavel~Zalesskii}{University of Brasilia}}
     {Zalesskii: Profinite topology on arithmetic groups}
  Let $G$ be a group. We can make  $G$ into a topological group by
  considering all normal finite index subgroups of $G$ as a
  fundamental system of neighborhoods of the identity. This topology
  is called the {\it profinite topology} on $G$.

  \begin{question}\label{Zal-qu-1}
    How strong is the profinite topology?
  \end{question}
  \begin{question}\label{Zal-qu-2}
    To what extent does the profinite completion
    \[
      \widehat G
      =
      \lim\limits_{\displaystyle\longleftarrow}{}_{N\triangleleft_f G}\  G/N
    \]
    determine $G$?
  \end{question}

  \begin{definition}\label{Zal-def-1}
    If every finitely generated subgroup of $G$ is closed in the profinite
    topology, then $G$ is called \notion{subgroup separable}.
  \end{definition}
  \begin{definition}\label{Zal-def-2}
    If the conjugacy class of every element is closed, then $G$ is
    called \notion{conjugacy separable}.
  \end{definition}

  \paragraph{Cohomological aspect of Question~\ref{Zal-qu-2}.} According to
  J.-P. Serre \cite{Zal-Serre}, a group $G$ is called \notion{good} if
  $G\longrightarrow \widehat G$ induces isomorphisms $H^n(\widehat
  G,M)\longrightarrow H^n(G,M)$ for every finite $G$-module $M$.

  Subgroup separability, conjugacy separability, goodness are
  indications (or features) of \emph{strong} profinite topology.

  The term ``strong profinite topology'' has a precise meaning for
  $S$-arithmetic groups. Namely, the profinite topology on an
  $S$-arithmetic group $\Gamma$ is \emph{strong} if $\Gamma$ does not have
  the Congruence Subgroup Property.

  \begin{remark}
    If an $S$-arithmetic group $\Gamma$ has {\small CSP} then it
    is \emph{not} subgroup separable, conjugacy separable or
    good.
  \end{remark}

  \begin{conjecture}
    If an $S$-arithmetic group $\Gamma$ does not have {\small CSP}
    then $\Gamma$ is conjugacy separable and subgroup separable.
    If in addition the characteristic of the ground field is
    zero, then $\Gamma$ is good.
  \end{conjecture}

  \paragraph{Supporting result.} The conjecture is true for
  arithmetic lattices in $SL_2(\CCC)$.

  Recent progress in the study of 3-manifolds allows one to deduce from results of
  Wilton-Zalesskii \cite{Zal-wilton_profinite_2010} that the
  fundamental group of a compact 3-manifold is good and from results
  of Hamilton-Wilton-Zalesskii \cite{Zal-HWZ} that it is conjugacy
  separable.
\endgroup

\newpage\addtocontents{toc}{\protect\bigskip}
\section*{References}
 	\markright{\textsc{References}}
   	\thispagestyle{plain}

\end{document}